\documentclass[12pt]{amsart}
\usepackage{amssymb,graphicx}
\usepackage{amsfonts}
\usepackage{latexsym}
\usepackage{amscd}
\usepackage{xypic}
\usepackage[all]{xy}
\newtheorem{prop}{Proposition}[section]
\newtheorem{lem}[prop]{Lemma}
\newtheorem{thm}[prop]{Theorem}
\newtheorem{cor}[prop]{Corollary}
\newtheorem{definition}[prop]{Definition}

\theoremstyle{definition}

\theoremstyle{remark}
\newtheorem{rem}[prop]{Remark}
\newtheorem*{question}{Question}

\numberwithin{equation}{section}

\newcommand{\C}{{\mathbb C}}
\newcommand\eqdef{{\;\overset{\mbox{\scriptsize def}}{=}\;}}
\newcommand{\Cs}{{$C^*$-al\-ge\-bra}}
\newcommand{\sh}{{$^*$-ho\-mo\-mor\-phism}}
\newcommand{\cP}{{\mathcal P}}

\begin{document}

\title{AF-embeddings into C*-algebras of real rank zero}

\author{Francesc Perera \& Mikael R\o rdam}
\address{Department of Pure Mathematics, Queen's University
Belfast, Belfast BT7 1NN, Northern Ireland}
\email{perera@qub.ac.uk}
\address{\emph{Current address}: Departament de Matem\`atiques, Universitat Aut\`onoma de Barcelona, 08193 Bellaterra, Barcelona, Spain}
\email{perera@mat.uab.es}
\address{Department of Mathematics, University of Southern Denmark,
  Campusvej 55, 5230 Odense M, Denmark} \email{mikael@imada.sdu.dk}

\date{October 2003}

\thanks{Research supported by the Danish Research Council, SNF, and the
Nuffield Foundation}


\maketitle

\begin{abstract}It is proved that every separable $C^*$-algebra of
  real rank zero contains an AF-sub-$C^*$-algebra such that the
  inclusion mapping induces an isomorphism of the ideal
  lattices of the two $C^*$-algebras and such that every projection in
  a matrix algebra over the large $C^*$-algebra is equivalent to a
  projection in a matrix 
  algebra over the AF-sub-$C^*$-algebra. This result is proved at the
  level of monoids, using that the monoid of Murray-von Neumann
  equivalence classes of projections in a $C^*$-algebra of real rank
  zero has the refinement property. As an application of our result,
  we show that given a unital
  $C^*$-algebra $A$ of real rank zero and a natural number $n$, then there is a unital
  $^*$-homomorphism $M_{n_1} \oplus \cdots \oplus M_{n_r} \to A$ for
  some natural numbers $r,n_1, \dots,n_r$ with $n_j \ge n$ for all $j$
  if and only if $A$ has no representation of dimension less than $n$.
\end{abstract}

\section{Introduction}
\noindent
Many properties of $C^*$-algebras, that can be
traced back to properties of high dimensional spaces, have recently
been found in simple $C^*$-algebras. These properties include
having arbitrary stable rank (Villadsen, \cite{Vil:sr=n}); perforation of the
ordered $K_0$ group (Villadsen, \cite{Vil:perforation}); and
existence of a finite projection that becomes infinite if doubled (the
last named author, \cite{Ror:simple}). Examples of simple nuclear
$C^*$-algebras with these properties appear to be in disagreement with
the classification conjecture by Elliott.

In contrast to the situation of (arbitrary) simple $C^*$-algebras,
no such exotic phenomenae have so far been found among
$C^*$-algebras of real rank zero (first studied by Brown and
Pedersen in \cite{bp}). At this juncture we do not know if this
lack of examples is because no such examples exist of if it is
because new constructions are needed to create them. With this
paper we wish to add to the knowledge of $C^*$-algebras of real
rank zero in the hope to shed some new light on their structure.

The AF-embedding theorem (explained in the abstract) extends a result
by H.\ Lin, \cite{lin}. He shows that if $A$ is a $C^*$-algebra of
real rank zero
and stable rank one and if $K_0(A)$ is a dimension group, then $A$
contains an AF-sub-$C^*$-algebra $B$ such
that the inclusion mapping induces an isomorphism $K_0(B) \to K_0(A)$
that maps the dimension range of $B$ onto that of $A$. (Such an
inclusion will necessarily induce an isomorphism between the ideal
lattices of $A$ and $B$.) In our result we do not require $A$ to have
stable rank one, nor do we require $K_0$ to be a dimension group.

The AF-embedding result shows in particular that the primitive
ideal space of an arbitrary separable $C^*$-algebra of real rank zero is
homeomorphic to that of the embedded AF-algebra. By a theorem of
Kirchberg in which $\mathcal{O}_2$-absorbing $C^*$-algebras are
classified by their primitive ideal spaces one can conclude that
the tensor product of an arbitrary separable, nuclear stable
$C^*$-algebra by $\mathcal{O}_2$ is isomorphic to the tensor
product of an AF-algebra by $\mathcal{O}_2$
(Corollary~\ref{cor:O_2}). It also follows from an older theorem
by Bratteli and Elliott that the primitive ideal space of a
separable real rank zero $C^*$-algebra is realized by an
AF-algebra (see Theorem~\ref{brat-ell}).

In Section~5 we use the AF-embedding result to derive some
divisibility properties of real rank zero $C^*$-algebras (one of which
is explained in the abstract). We say that a $C^*$-algebra $A$ is
weakly divisible if for every projection $p$ in $A$ there is a
unital $^*$-homomorphism $M_{n_1} \oplus \cdots \oplus M_{n_r} \to
pAp$ whenever $\gcd\{n_1, \dots, n_r\} =1$. It is shown in
Theorem~\ref{thm:divisibility} that a separable real rank
zero algebra is weakly divisible if and only if no representation
of $A$ has non-zero intersection with the compact operators.

The proof of the main result has two parts. In the first part,
treated in Section~3, one solves the problem at the algebraic
level of monoids. In more detail, one associates to each
$C^*$-algebra $A$ a monoid $V(A)$, consisting of Murray-von
Neumann equivalence classes of projections in $M_\infty(A)$ (the
union of all matrix algebras over $A$). The monoid $V(A)$ has the
so-called refinement property (see Section~3) when $A$ is of real
rank zero. It is shown in Section~3 that if $M$ is a monoid with
the refinement property then there is a dimension monoid $\Delta$
and a surjective monoid morphism $\alpha \colon \Delta \to M$ such
that for all $x,y \in \Delta$ if $\alpha(x) \le m\alpha (y)$ for
some natural number $m$, then $x \le ny$ for some natural number
$n$. (A morphism $\alpha$ with this property will induce an
isomorphism of the ideal lattices of the two monoids $\Delta$ and
$M$, see below for the precise definitions.) This result also allows us
to recover a representation result of distributive lattices due to
Goodearl and Wehrung (see~\cite{gw}).

The second part of the proof of the embedding result, treated in
Section~4, is fairly standard and consists of lifting a monoid
morphism $V(B) \to V(A)$, where $B$ is an AF-algebra and $A$ is an
arbitrary (stable) $C^*$-algebra to a $^*$-homomorphism. (Note
that we do not require $A$ to have the cancellation property.)

\section{Isomorphism of ideal lattices}
\noindent
In this short section we will establish a basic result that
motivates the purely monoid theoretic work carried out in the
sequel.


For any $C^*$-algebra $A$, we will in this paper use the notation
$L(A)$ to refer to the lattice of two-sided closed ideals of $A$.
Given a $^*$-homomorphism $\varphi\colon A\to B$ between
$C^*$-algebras $A$ and $B$, there is a natural way of relating the
ideal lattice of $A$ with that of $B$. This is given by the set
map $\varphi^{-1}\colon L(B)\to L(A)$, defined by taking
pre-images of ideals.

We shall be interested in the situation when the map
$\varphi^{-1}$ is an isomorphism, and we refer to that by saying
that $\varphi$ \emph{induces an isomorphism of ideal lattices}.
Note that such a $^*$-homomorphism $\varphi$ will always be
injective (as the zero ideal is contained in its kernel).

For any $C^*$-algebra $A$, let us denote by $\mathrm{Proj}(A)$ the
set of its projections. If $p$, $q\in \mathrm{Proj}(A)$, we denote
as usual $p\sim q$ to mean that they are \emph{Murray-von Neumann
equivalent}, that is, $p=vv^*$ and $q=v^*v$ for a (partial
isometry) $v$ in $A$. Let $M_{\infty}(A)=\varinjlim M_n(A)$, under
the mappings $M_n(A)\to M_{n+1}(A)$, $x\mapsto \left(\begin{smallmatrix} x & 0 \\
0 & 0
\end{smallmatrix}\right)$. Denote by $V(A)$ the set of Murray-von
Neumann equivalence classes $[p]$ of projections $p$ coming from
$M_{\infty}(A)$. This is an abelian monoid that admits a natural
preordering, namely $x\leq y$ if there is a $z$ in $V(A)$ such that
$x+z=y$. If $A$ is unital, then the Grothendieck group of $V(A)$
is $K_0(A)$. Any $^*$-homomorphism $\varphi\colon A\to B$ induces a
monoid morphism $V(\varphi)\colon V(A)\to V(B)$ by
$V(\varphi)([p])=[\varphi(p)]$.

Recall that a unital $C^*$-algebra $A$ has \emph{real rank zero}
provided that every self-adjoint element can be approximated
arbitrarily well by self-adjoint, invertible elements. If $A$ is
non-unital, then $A$ is said to have real rank zero if its minimal
unitization $\widetilde A$ has real rank zero (see \cite{bp}).
Hereditary subalgebras of $C^*$-algebras of real rank zero (in
particular ideals) have also real rank zero by \cite[Corollary
2.8]{bp}, and hence they can be written as the closed linear span
of their projections. It follows from this that two closed ideals
in a $C^*$-algebra of real rank zero will be equal precisely when
they have the same projections.

Now let $M$ be an abelian monoid, equipped with its natural
(algebraic) ordering. We say that a submonoid $I$ of $M$ is an
\emph{order-ideal} (or \emph{o-ideal}) if $x+y\in I$ if and
only if $x$, $y\in I$. (Equivalently, if $x\leq y$, and $y\in I$,
then $x\in I$.) Denote the ideal lattice of a monoid $M$ by
$L(M)$. If $f\colon M\to N$ is a monoid morphism, then we say that
$f$ \emph{induces an isomorphism of ideal lattices} if the set map
$f^{-1}\colon L(N)\to L(M)$ is an isomorphism. In contrast to the
$C^*$-algebra case, a morphism that induces an isomorphism of
ideal lattices need not be injective, due to the absence of
kernels in the category of abelian monoids. Let $M$ be an abelian monoid, and
let $x$, $y\in M$. In order to ease the notation in this and the next
sections, we shall denote $x\propto y$ to mean that $x\leq ny$ for
some $n$ in $\mathbb N$,
or, equivalently, $x$ belongs to the order-ideal in $M$ generated by $y$.

\begin{lem}
\label{idealmonoid} Let $M$ and $N$ be abelian monoids and
let $\alpha\colon M\to N$ be a surjective monoid morphism. Suppose
that $x\propto y$ for any $x$, $y$ in $M$ for which $\alpha(x)\propto
\alpha(y)$. Then $\alpha$ induces an isomorphism of ideal
lattices.
\end{lem}

\begin{proof}
Assume that, for order-ideals $I$ and $J$ of $N$, we have
$\alpha^{-1}(I)=\alpha^{-1}(J)$. Let $x\in I$. Then there is $y$
in $M$ such that $\alpha(y)=x$, and so $y\in
\alpha^{-1}(I)=\alpha^{-1}(J)$. This says $x=\alpha(y)\in J$. Thus
$J\subseteq I$ and by symmetry $J=I$. This proves that
$\alpha^{-1}\colon L(N)\to L(M)$ is injective.

Next, if $I$ is an order-ideal of $M$, let $J=\{y\in N\mid y\propto
\alpha(x)\text{ for some }x\text{ in }I\}$, which is an
order-ideal in $N$. Clearly $I\subseteq \alpha^{-1}(J)$. If $x\in
\alpha^{-1}(J)$, then $\alpha(x)\in J$, so $\alpha(x)\propto
\alpha(y)$ for some $y$ in $I$. Our assumption on $\alpha$ implies that
$x\propto y$, hence $x\in I$.
\end{proof}

\noindent Given a $C^*$-algebra $A$, there is a natural mapping
$L(A)\to L(V(A))$, given by the rule $I\mapsto V(I)$. If $A$ has
real rank zero, then this mapping is a lattice isomorphism, as
proved by Zhang in~\cite{zh}. We shall use this fact in the
following lemma.
\begin{lem}
\label{isomonoid} Let $A$ and $B$ be $C^*$-algebras with real rank
zero, and let $\varphi\colon A\to B$ be a $^*$-homomorphism. Then
$\varphi$ induces an isomorphism of ideal lattices if and only if
$V(\varphi)$ does.
\end{lem}

\begin{proof}
First we note that, if $I$ is an ideal of $B$, then
\[
V(\varphi)^{-1}(V(I))=V(\varphi^{-1}(I))\,.
\]
Indeed, let $x\in V(\varphi)^{-1}(V(I))$. Write $x=[p]$, where
$p\in M_n(A)$, for some $n$. Then $V(\varphi)(x)=[\varphi(p)]\in
V(I)$, so that $\varphi(p)\in M_n(I)$. This means that
$p\in\varphi^{-1}(M_n(I))$, hence $x=[p]\in V(\varphi^{-1}(I))$.
The converse inclusion is similar.

This says that the diagram
\[\begin{CD} L(A) @>{\cong}>> L(V(A))\\
@A{\varphi^{-1}}AA  @AA{V(\varphi)^{-1}}A\\
L(B) @>{\cong}>> L(V(B))
\end{CD}\]
is commutative. Therefore the conclusion of the lemma follows at once.
%
%
\end{proof}

\noindent Let us summarize our observations in the following:
\begin{cor}
\label{ideal2} Let $A$ and $B$ be $C^*$-algebras with real rank
zero, and let $\varphi\colon A\to B$ be a $^*$-homomorphism.
Suppose that $V(\varphi)\colon V(A)\to V(B)$ is onto and that
$x\propto y$ for any $x$, $y$ in $V(A)$ for which $V(\varphi)(x)\propto
V(\varphi)(y)$. Then $\varphi$ is injective and induces an
isomorphism of ideal lattices.
\end{cor}

\section{Mapping dimension monoids onto refinement monoids}

In order to use Corollary~\ref{ideal2} effectively, we need to
gain insight into surjective mappings between monoids that arise
as equivalence classes of projections of $C^*$-algebras with real
rank zero. As indicated in the introduction, one of them will
actually come from an $AF$-algebra, and hence will be a dimension
monoid (see below).

An abelian monoid $M$ is said to be \emph{conical} provided
$x+y=0$ precisely when $x=y=0$. Note that if $A$ is a
$C^*$-algebra, then $V(A)$ is a conical monoid. Therefore, in this
section, all monoids are assumed to be conical.

Given a natural number $r$, we shall refer to the monoid $(\mathbb
Z^+)^r$ as a \emph{simplicial} monoid (of rank $r$). Observe that
such a monoid has a \emph{canonical basis}, namely the one
obtained by setting $e_i=(0,\ldots,0,1,0,\ldots,0)$ (with $1$ in
the $i$-th position), for $1\leq i\leq r$. These basis elements
are precisely the minimal non-zero elements of $(\mathbb Z^+)^r$.
A monoid that can be written as an inductive limit of a sequence
of (finitely generated) simplicial monoids and monoid morphisms
will be called a \emph{dimension monoid}. (The reason for this
terminology is that dimension monoids are precisely the
positive cones of dimension groups.)

The structure of $V(A)$ for a general $C^*$-algebra $A$ with real
rank zero can be more intricate than just being a dimension
monoid. Still, any such $V(A)$ will enjoy the following important
property:

An abelian monoid $M$ is termed a \emph{refinement} monoid (see,
e.g. \cite{w} or \cite{dob}) if whenever $x_1+x_2=y_1+y_2$ in $M$,
then there exist elements $z_{ij}$ in $M$ such that
$x_i=z_{i1}+z_{i2}$ and $y_i=z_{1i}+z_{2i}$, for $i=1$, $2$. That
$V(A)$ is a refinement monoid for any $C^*$-algebra $A$ of real
rank zero is proved in \cite[Lemma 2.3]{ap}, based on work of
Zhang (\cite{zh}).

The next few lemmas are of a technical nature and will be used to
assemble the proof of the main monoid-theoretical result
(Theorem~\ref{mainmonoid}). Applications to operator algebras will
be given in the next sections.

Let $\Delta$ be a simplicial monoid of rank $r$ and canonical
basis $e_1,\ldots,e_r$. If $I$ is a finite subset of
$\{1,\ldots,r\}$, we denote $e_I=\sum_{i\in I}e_i$. If
$I=\{1,\ldots,r\}$, then we shall denote $e_{\Delta}=e_I$, and if
$I=\emptyset$, then $e_I=0$.

In the next lemma we shall make a critical use of the refinement
property. One key step in the proof is provided by a result of
F.~Wehrung in \cite{w}. Recall that if $M$ is an abelian monoid
and $x$, $y\in M$, we use the notation $x\propto y$ to mean $x\leq ny$
for some $n$ in $\mathbb N$.
\begin{lem}
\label{monoid1} Let $\Delta$ be a simplicial monoid of rank $r$.
Let $M$ be a refinement monoid and $\alpha\colon \Delta\to M$ be a
monoid morphism. Suppose that
\[
\alpha(e_1)\propto \alpha(e_2+\cdots+e_r)
\]
in $M$. Then there exist a simplicial monoid $\Delta'$ (with
$\mathrm{rank}(\Delta')\geq 2(r-1)$), monoid morphisms
$\beta\colon\Delta\to\Delta'$ and $\alpha'\colon\Delta'\to M$ such
that $\alpha(\Delta)\subseteq\alpha'(\Delta')$, the diagram
\[\begin{CD} \Delta @>{\alpha}>> \alpha(\Delta)\\
@V{\beta}VV  @VV{\iota}V\\
\Delta' @>{\alpha'}>> \alpha'(\Delta')
\end{CD}\]
is commutative, and
$\beta(e_1)\propto\beta(e_2+\cdots+e_r)=e_{\Delta'}$ in $\Delta'$.
\end{lem}

\begin{proof}
Our assumption means that $\alpha(e_1)\leq
n\alpha(e_2+\cdots+e_r)$ for some $n$ in $\mathbb N$.

Since $M$ is a refinement monoid, we can use \cite[Lemma 1.9]{w}
in order to find elements $y_0,\ldots, y_n$ in $M$ such that
\[
\alpha(e_1)=\sum_{j=1}^n jy_j\,,\mbox{ and }
\alpha(e_2+\cdots+e_r)=\sum_{j=0}^n y_j.
\]
Another use of the refinement property (applied to the second
identity above) yields elements $x_{ij}$ in $M$, for $2\leq i\leq
r$ and $0\leq j \leq n$ such that
\[
\alpha(e_i)=\sum_{j=0}^n x_{ij}\,,\mbox{ for } i=2,\ldots,r
\]
and
\[
y_j=\sum_{i=2}^r x_{ij}\,,\mbox{ for } j=0,\ldots,n\,.
\]
Next, let $\Delta'$ be the simplicial monoid of rank $(n+1)(r-1)$,
and denote its canonical basis by $(e_{ij})$, where $2\leq i\leq
r$ and $0\leq j\leq n$. Now we can define $\alpha'\colon\Delta'\to
M$ by $\alpha'(e_{ij})=x_{ij}$. Define
$\beta\colon\Delta\to\Delta'$ by
\[
\beta(e_1)=\sum_{i=2}^r\sum_{j=1}^n je_{ij}\,,\mbox{ and
}\beta(e_i)=\sum_{j=0}^n e_{ij}\,\mbox{ for }2\leq i\leq r\,.
\]
Then
\[
\alpha'(\beta(e_1))=\sum_{i=2}^r\sum_{j=1}^n
jx_{ij}=\sum_{j=1}^njy_j=\alpha(e_1)
\]
and
\[
\alpha'(\beta(e_i))=\sum_{j=0}^n x_{ij}=\alpha(e_i)\,.
\]
Thus $\alpha'\circ\beta=\alpha$, and in particular we see that
$\alpha(\Delta)=(\alpha'\circ\beta)(\Delta)\subseteq\alpha'(\Delta')$.
Finally,
\[
\beta(e_1)=\sum_{i=2}^r\sum_{j=0}^n je_{ij}\leq
n\sum_{i=2}^r\sum_{j=0}^n e_{ij}=n\beta(e_2+\cdots+e_r)\,,
\]
so $\beta(e_1)\propto\beta(e_2+\cdots+e_r)$. Also
$\beta(e_2+\cdots+e_r)=\sum_{i=2}^r\sum_{j=0}^n
e_{ij}=e_{\Delta'}$.
\end{proof}

\noindent Of course, the sets of indices $\{1\}$ and
$\{2,\ldots,r\}$ in the above lemma can be replaced by $\{j\}$ and
$\{1,\ldots,r\}\setminus\{j\}$ to achieve a similiar conclusion.

The following easy fact will be used tacitly a number of times in
what follows.
\begin{rem}
Let $\Delta$ be a simplicial monoid of rank $r$ and canonical
basis $\{e_j\}_{1\leq j\leq r}$, and let $\alpha\colon\Delta\to M$
be a monoid morphism, where $M$ is an abelian monoid. If
$I\subseteq\{1,\ldots,t\}$ and $y\in M$, then $\alpha(e_I)\propto y$
if and only if $\alpha(e_i)\propto y$ for all $i$ in $I$. Indeed, to
prove the less trivial part, assume that there are natural numbers
$n_i$, for $i$ in $I$, such that $\alpha(e_i)\leq n_iy$ for all
$i$. Then $\alpha(e_I)=\sum_{i\in I}\alpha(e_i)\leq (\sum_i
n_i)y\propto y$.
\end{rem}

\begin{lem}
\label{inductivestep} Let $\Delta$ be a simplicial monoid of rank
$r$. Let $I$, $J$ be non-empty subsets of $\{1,\ldots,r\}$ such
that $I\cap J=\emptyset$. Suppose that there is a refinement
monoid $M$ and a monoid morphism $\alpha\colon\Delta\to M$ such
that $\alpha(e_J)\propto\alpha(e_I)$.

Let $j\in J$. Then there are a simplicial monoid $\Delta'$ with
$\mathrm{rank}(\Delta')=s\geq 2|I|+|J|-1$ (and canonical basis
denoted by $\{f_k\}_{k=1}^s$), monoid morphisms
$\alpha'\colon\Delta'\to M$ and $\beta\colon\Delta\to\Delta'$ such
that the following conditions hold:
\begin{enumerate}\itemsep=2mm
\item[(i)] $\alpha(\Delta)\subseteq\alpha'(\Delta')$\,.

\item[(ii)] The diagram
\[\begin{CD} \Delta @>{\alpha}>> \alpha(\Delta)\\
@V{\beta}VV  @VV{\iota}V\\
\Delta' @>{\alpha'}>> \alpha'(\Delta')
\end{CD}\]
is commutative.

\item[(iii)] $\beta(e_k)=f_k$ for all $k\in J\setminus\{j\}$\,.

\item[(iv)] There is a subset $K$ of $\{1,\ldots,s\}$ such that
$K\cap(J\setminus\{j\})=\emptyset$, $\beta(e_j)\propto\beta(e_I)=f_K$,
and $\alpha'(f_{J\setminus\{j\}})\propto\alpha'(f_K)$\,.
\end{enumerate}
\end{lem}

\begin{proof}
Write $\Delta=\Delta_1\oplus\Delta_2$ (i.e.
$\Delta=\Delta_1+\Delta_2$ and $\Delta_1\cap\Delta_2=0$). Here
$\Delta_1$ is the simplicial submonoid of $\Delta$ with basis
elements indexed by $\{1,\ldots,r\}\setminus (I\cup\{j\})$, so it
has rank at least $|J|-1$, and $\Delta_2$ is spanned by the
remaining basis elements. Define $\alpha_2\colon\Delta_2\to M$ by
restriction of $\alpha$. By assumption, we have that
$\alpha_2(e_j)\propto\alpha_2(e_I)$.

By Lemma~\ref{monoid1}, there are a simplicial monoid $\Delta_2'$
(of rank at least $2|I|$), monoid morphisms
$\alpha_2'\colon\Delta_2'\to M$,
$\beta_2\colon\Delta_2\to\Delta_2'$ such that
$\alpha_2(\Delta_2)\subseteq\alpha_2'(\Delta_2')$, the diagram
\[\begin{CD} \Delta_2 @>{\alpha_2}>> \alpha_2(\Delta_2)\\
@V{\beta_2}VV  @VV{\iota}V\\
\Delta_2' @>{\alpha_2'}>> \alpha_2'(\Delta_2')
\end{CD}\]
is commutative, and $\beta_2(e_j)\propto\beta_2(e_I)$.

Set $\Delta'=\Delta_1\oplus\Delta_2'$ (this is an external direct
sum). The canonical basis $\{f_k\}$ of $\Delta'$ is obtained by
taking the union of the canonical bases of $\Delta_1$ and
$\Delta_2'$.

Note that
$\mathrm{rank}(\Delta')=\mathrm{rank}(\Delta_1)+\mathrm{rank}(\Delta_2')\geq
2|I|+|J|-1$. Define $\beta\colon\Delta\to\Delta'$ by
$\beta_{|\Delta_1}=\mathrm{id}$, and by
$\beta_{|\Delta_2}=\beta_2$. Define $\alpha'\colon\Delta'\to M$ by
$\alpha'_{|\Delta_1}=\alpha$ and $\alpha'_{|\Delta_2'}=\alpha_2'$.
Observe that we get
\[
\alpha(\Delta)=\alpha(\Delta_1)+\alpha(\Delta_2)=\alpha(\Delta_1)+\alpha_2(\Delta_2)\subseteq
\alpha(\Delta_1)+\alpha_2'(\Delta_2')=\alpha'(\Delta_1)+\alpha'(\Delta_2')=\alpha'(\Delta')\,,
\]
and also the following commutative diagram:
\[\begin{CD} \Delta  &=\Delta_1\oplus \Delta_2 @>{\alpha}>> \alpha(\Delta)\\
 &  @V{\beta}VV  @VV{\iota}V\\
\Delta' & =\Delta_1\oplus\Delta_2' @>{\alpha'}>> \alpha'(\Delta')
\end{CD}\]
Hence conditions (i) and (ii) are fulfilled. By construction
$\beta(e_k)=f_k$ for any $k$ in $J\setminus\{j\}$ (as they belong
to $\Delta_1$). Thus condition (iii) also holds.

Note that $\beta(e_j)=\beta_2(e_j)\propto\beta_2(e_I)=\beta(e_I)$.

Let $K$ be the set of indices corresponding to a basis of
$\Delta_2'$. Clearly, $K\cap(J\setminus\{j\})=\emptyset$. By
definition of $\beta$ and Lemma~\ref{monoid1}, we see that
$\beta(e_I)=\beta_2(e_I)=f_{\Delta_2'}=f_{K}$.

Finally,
$\alpha'(f_{J\setminus\{j\}})=\alpha'(\beta((e_{J\setminus\{j\}})))=\alpha(e_{J\setminus\{j\}})\propto\alpha(e_I)=\alpha'(\beta(e_I))=\alpha'(f_K)$.
This verifies condition (iv).
\end{proof}

\begin{lem}
\label{induct2}Let $\Delta$ be a simplicial monoid of rank $r$.
Let $I$, $J$ be non-empty subsets of $\{1,\ldots, r\}$ such that
$I\cap J=\emptyset$. Suppose there is a refinement monoid $M$ and
a monoid morphism $\alpha\colon\Delta\to M$ such that
$\alpha(e_J)\propto\alpha(e_I)$.

Then there are a simplicial monoid $\Delta'$, monoid morphisms
$\alpha'\colon\Delta'\to M$ and $\beta\colon\Delta\to\Delta'$ such
that $\alpha(\Delta)\subseteq\alpha'(\Delta')$, the diagram
\[\begin{CD} \Delta @>{\alpha}>> \alpha(\Delta)\\
@V{\beta}VV  @VV{\iota}V\\
\Delta' @>{\alpha'}>> \alpha'(\Delta')
\end{CD}\]
is commutative, and $\beta(e_J)\propto\beta(e_I)$.
\end{lem}

\begin{proof}
We proceed by induction on $|J|$. If $|J|=1$, the conclusion is
provided by (the argument of) Lemma~\ref{monoid1}.

Assume that $|J|\geq 2$ and that the conclusion holds for
cardinals smaller than $|J|$. Let $j$ be the first element in $J$
(in the natural ordering). By Lemma~\ref{inductivestep}, there are
a simplicial monoid $\Delta_1'$ with $\mathrm{rank}(\Delta')=s\geq
2|I|+|J|-1$ (and canonical basis denoted by $\{f_k\}_{k=1}^s$),
monoid morphisms $\alpha_1'\colon\Delta_1'\to M$ and
$\beta_1\colon\Delta\to\Delta_1'$ such that
$\alpha(\Delta)\subseteq\alpha_1'(\Delta')$\,, the diagram
\[\begin{CD} \Delta @>{\alpha}>> \alpha(\Delta)\\
@V{\beta_1}VV  @VV{\iota}V\\
\Delta_1' @>{\alpha_1'}>> \alpha_1'(\Delta')
\end{CD}\]
is commutative, $\beta_1(e_k)=f_k$ for all $k\in
J\setminus\{j\}$\,, $\beta_1(e_j)\propto\beta_1(e_I)=f_K$, and
$\alpha_1'(f_{J\setminus\{j\}})\propto\alpha_1'(f_K)$, for some subset
$K$ of $\{1,\ldots,s\}$ such that
$K\cap(J\setminus\{j\})=\emptyset$\,.

By the induction hypothesis applied to $\Delta_1'$, $K$,
$J\setminus\{j\}$, $\alpha_1'$, there are a simplicial monoid
$\Delta'$, monoid morphisms $\alpha'\colon\Delta'\to M$ and
$\beta_2\colon\Delta_1'\to\Delta'$ such that
$\alpha_1'(\Delta_1')\subseteq\alpha'(\Delta')$, the diagram
\[\begin{CD} \Delta_1' @>{\alpha_1'}>> \alpha_1'(\Delta_1')\\
@V{\beta_2}VV  @VV{\iota}V\\
\Delta' @>{\alpha'}>> \alpha'(\Delta')
\end{CD}\]
is commutative, and $\beta_2(f_{J\setminus\{j\}})\propto\beta_2(f_K)$.

Next, define $\beta=\beta_2\circ\beta_1$. We have that
\begin{align*}
\beta(e_J) &
=\beta_2(\beta_1(e_J))=\beta_2(\beta_1(e_j))+\beta_2(\beta_1(e_{J\setminus\{j\}}))
\\
&
\propto\beta_2(\beta_1(e_I))+\beta_2(f_{J\setminus\{j\}})\\
&
\propto\beta_2(\beta_1(e_I))+\beta_2(f_K)\\
&
=\beta_2(\beta_1(e_I))+\beta_2(\beta_1(e_I))=2\beta_2(\beta_1(e_I))\propto\beta(e_I)\,,
\end{align*}
as desired.
\end{proof}

\noindent Let $\Delta$ be a simplicial monoid of rank $r$ with
canonical basis denoted by $\{e_i\}_{1\leq i\leq r}$\,. For any
$j$, we have the coordinate monoid morphism
$\pi_j\colon\Delta\to\mathbb Z^+$, i.e. $\pi_j(x)=\alpha_j$
provided that $x=\sum_{i=1}^r\alpha_i e_i$ (with all $\alpha_i$ in
$\mathbb Z^+$).

For any element $x$ in $\Delta$, we define its \emph{support} as
the following subset of $\{1,\ldots,r\}$
\[
\mathrm{supp} (x)=\{j\mid \pi_j(x)\neq 0\}\,.
\]
\begin{lem}
\label{supports} Let $\Delta$ be a simplicial monoid of rank $r$
and canonical basis $\{e_j\}_{1\leq j\leq r}$. Suppose that $M$ is
an abelian monoid and that $\alpha\colon\Delta\to M$ is a monoid
morphism. For elements $x$ and $y$ in $\Delta$ with
$J=\mathrm{supp}(x)$ and $I=\mathrm{supp}(y)$, the following
conditions are equivalent:
\begin{enumerate}\itemsep=2mm
\item[(i)] $\alpha(x)\propto\alpha(y)$\,.

\item[(ii)] $\alpha(e_J)\propto\alpha(e_I)$\,.

\item[(iii)] $\alpha(e_{J\setminus I})\propto\alpha(e_I)$\,.
\end{enumerate}
\end{lem}

\begin{proof}
Write $x=\sum_{j\in J}\alpha_j e_j$ and $y=\sum_{i\in I}\beta_i
e_i$, where $\alpha_j\geq 1$ and $\beta_i\geq 1$ for all $j$ and
$i$.

(i) $\Rightarrow$ (ii). Assume $\alpha(x)\propto\alpha(y)$. Then, for
any $j$ in $J$ we have that
\[
\alpha(e_j)\leq\alpha_j\alpha(e_j)\leq\alpha(x)\propto\alpha(y)=\sum_{i\in
I}\beta_i\alpha(e_i)\leq\max\{\beta_i\mid i\in
I\}\alpha(e_I)\propto\alpha(e_I)\,.
\]
Hence $\alpha(e_J)\propto\alpha(e_I)$.

(ii) $\Rightarrow$ (iii). Trivial.

(iii) $\Rightarrow$ (i). Assume that $\alpha(e_{J\setminus
I})\propto\alpha(e_I)$. Then
\begin{align*}
\alpha(x)=\sum_{j\in J}\alpha_j\alpha(e_j) & =\sum_{j\in
J\setminus I}\alpha_j\alpha(e_j)+\sum_{j\in J\cap I}\alpha_j
\alpha(e_j)\\
& \leq\max\{\alpha_j\mid j\in J\setminus I\}\alpha(e_{J\setminus
I})+\max\{\alpha_j\mid j\in J\cap I\}\alpha(e_{J\cap I})\\
& \propto 2\alpha(e_I)\propto\sum_{i\in I}\beta_i\alpha(e_i)=\alpha(y)\,.
\end{align*}
\end{proof}

\begin{lem}
\label{arbitrary elements} Let $\Delta$ be a simplicial monoid of
rank $r$. Let $M$ be a refinement monoid and
$\alpha\colon\Delta\to M$ be a monoid morphism. Suppose that $x$
and $y$ in $\Delta$ satisfy $\alpha(x)\propto\alpha(y)$.

Then there are a simplicial monoid $\Delta'$, monoid morphisms
$\alpha'\colon\Delta'\to M$ and $\beta\colon\Delta\to\Delta'$ such
that $\alpha(\Delta)\subseteq\alpha'(\Delta')$, the diagram
\[\begin{CD} \Delta @>{\alpha}>> \alpha(\Delta)\\
@V{\beta}VV  @VV{\iota}V\\
\Delta' @>{\alpha'}>> \alpha'(\Delta')
\end{CD}\]
is commutative, and $\beta(x)\propto\beta(y)$.
\end{lem}

\begin{proof}
Write $J=\mathrm{supp}(x)$ and $I=\mathrm{supp}(y)$. We may
clearly assume that $J\neq\emptyset$. If $I=\emptyset$, then $y=0$
and our assumption means that $\alpha(x)\leq 0$. By conicality,
$\alpha(x)=0$. If $x=\sum_{j\in J}\alpha_j e_j\,$, then
$\alpha(e_j)=0$ for every $j$ in $J$. Take $\Delta'=\Delta$, and
set $\beta(e_j)=e_j$ if $j\notin J$, $\beta(e_j)=0$ otherwise.
Take $\alpha'=\alpha$. Clearly, $\alpha\circ\beta=\alpha$ and
$\beta(x)=0$.

So, we may assume that $I$ and $J$ are both non-empty. By
Lemma~\ref{supports}, our assumption means that
$\alpha(e_{J\setminus I})\propto\alpha(e_I)$. We can therefore use
Lemma~\ref{induct2} to find a simplicial monoid $\Delta'$, monoid
morphisms $\alpha'\colon\Delta'\to M$ and
$\beta\colon\Delta\to\Delta'$ such that
$\alpha(\Delta)\subseteq\alpha'(\Delta')$, the diagram
\[\begin{CD} \Delta @>{\alpha}>> \alpha(\Delta)\\
@V{\beta}VV  @VV{\iota}V\\
\Delta' @>{\alpha'}>> \alpha'(\Delta')
\end{CD}\]
is commutative, and $\beta(e_{J\setminus I})\propto\beta(e_I)$. Now
Lemma~\ref{supports} implies that $\beta (x)\propto\beta(y)$, as
desired.
\end{proof}

\begin{lem}
\label{solving} Let $\Delta$ be a simplicial monoid of rank $r$,
and let $\alpha\colon\Delta\to M$ be a monoid morphism, where $M$
is a refinement monoid. Then there are a simplicial monoid
$\Delta'$, monoid morphisms $\beta\colon\Delta\to\Delta'$ and
$\alpha'\colon\Delta'\to M$ such that
$\alpha(\Delta)\subseteq\alpha'(\Delta')$, the diagram
\[\begin{CD} \Delta @>{\alpha}>> \alpha(\Delta)\\
@V{\beta}VV  @VV{\iota}V\\
\Delta' @>{\alpha'}>> \alpha'(\Delta')
\end{CD}\]
is commutative, and such that $\beta(x)\propto\beta(y)$ whenever $x$
and $y$ are elements in $\Delta$ with $\alpha(x)\propto\alpha(y)$.
\end{lem}

\begin{proof}
Our Lemma~\ref{supports} reduces the comparison of $\alpha(x)$ and
$\alpha(y)$ (with respect to the relation $\propto$) to the comparison
of their respective supports. Therefore it suffices to arrange the
construction for a finite set of pairs written as $\{(x_i,y_i)\mid
1\leq i\leq n\}$ (for some $n$ in $\mathbb N$) such that
$\alpha(x_i)\propto\alpha(y_i)$ in $M$ for every $i$.

This is done by iteration of Lemma~\ref{arbitrary elements},
starting with the pair $(x_1,y_1)$. Notice that in this first
step, we get the following diagram
\[\begin{CD} \Delta @>{\alpha}>> \alpha(\Delta)\\
@V{\beta_1}VV  @VV{\iota}V\\
\Delta_1' @>{\alpha_1'}>> \alpha_1'(\Delta_1')
\end{CD}\]
For the next pair $(x_2,y_2)$, just observe that
$\alpha_1'(\beta_1(x_2))=\alpha(x_2)\propto\alpha(y_2)=\alpha_1'(\beta_1(x_2))$,
hence the process can be reiterated.
\end{proof}

\begin{prop}
\label{almostthere} Let $M=\{x_0,x_1,\ldots\}$ be a countable
refinement monoid. There is then a commutative diagram:
\begin{equation}
\label{diagram}
\begin{CD} \Delta_0 @>{\alpha_0}>> M_0\\
@V{\beta_0}VV  @VV{\iota}V\\
\Delta_1 @>{\alpha_1}>> M_1\\
@V{\beta_1}VV  @VV{\iota}V\\
\Delta_2 @>{\alpha_2}>> M_2\\
@V{\beta_2}VV  @VV{\iota}V\\
\vdots & & \vdots
\end{CD}
\end{equation}
such that $\{x_0,\ldots,x_j\}\subseteq M_j$ for all $j$, each
$\Delta_j$ is a simplicial monoid, the morphisms
$\alpha_j\colon\Delta_j\to M_j$ are surjective, and for each $j$
and each pair of elements $x$ and $y$ in $\Delta_j$ for which
$\alpha_j(x)\propto\alpha_j(y)$ in $M$, it follows that
$\beta_j(x)\propto\beta_j(y)$ in $\Delta_{j+1}$.
\end{prop}

\begin{proof}
Let $M_0=\langle x_0\rangle$ be the submonoid generated by $x_0$,
and let $\Delta_0$ be the simplicial monoid of rank $1$, with
basis $\{e_0\}$. Define $\alpha_0\colon\Delta_0\to M$ by
$\alpha_0(e_0)=x_0$. Clearly, $\alpha_0(\Delta_0)=M_0$.

We can now use Lemma~\ref{solving} to find a simplicial monoid
$\Delta_1'$, monoid morphisms $\beta_0'\colon\Delta_0\to\Delta_1'$
and $\alpha_1'\colon\Delta_1'\to M$ such that
$M_0\subseteq\alpha_1'(\Delta_1')$, the diagram
\[\begin{CD} \Delta_0 @>{\alpha_0}>> M_0\\
@V{\beta_0'}VV  @VV{\iota}V\\
\Delta_1' @>{\alpha_1'}>> \alpha_1'(\Delta_1')
\end{CD}\]
is commutative, and $\beta_0'(x)\propto\beta_0'(y)$ for any pair of
elements $x$, $y$ in $\Delta_0$ such that
$\alpha_0(x)\propto\alpha_0(y)$ (in $M$).

If $x_1\in\alpha_1'(\Delta_1')$, then set $\Delta_1=\Delta_1'$,
$\alpha_1=\alpha_1'$, $\beta_0=\beta_0'$ and
$M_1=\alpha_1(\Delta_1)$. Otherwise, set
$\Delta_1=\Delta_1'\oplus\langle e\rangle$, where $e$ is an
additional basis element, so that $\Delta_1$ is also a simplicial
monoid. Define $\alpha_1\colon\Delta_1\to M$ by
${\alpha_1}_{|\Delta_1'}=\alpha_1'$ and $\alpha_1(e)=x_1$, and set
$M_1=\alpha_1(\Delta_1)$, a submonoid of $M$. Define also
$\beta_0\colon\Delta_0\to\Delta_1$ as the composition of
$\beta_0'$ with the natural inclusion of $\Delta_1'$ in
$\Delta_1$.

Suppose that $\Delta_j$, $M_j$, $\alpha_j\colon\Delta_j\to M_j$
and $\beta_{j-1}\colon\Delta_{j-1}\to\Delta_j$ have been
constructed satisfying the requirements of our statement. Apply
Lemma~\ref{solving} to obtain $\Delta_{j+1}'$, morphisms
$\alpha_{j+1}'\colon\Delta_{j+1}'\to M$ and
$\beta_j'\colon\Delta_j\to\Delta_{j+1}'$ such that
$M_j\subseteq\Delta_{j+1}'$, the diagram
\[\begin{CD} \Delta_j @>{\alpha_j}>> M_j\\
@V{\beta_j'}VV  @VV{\iota}V\\
\Delta_{j+1}' @>{\alpha_{j+1}'}>> \alpha_{j+1}'(\Delta_{j+1}')
\end{CD}\]
is commutative, and $\beta_j'(x)\propto\beta_j'(y)$ for every pair of
elements $x$, $y$ in $\Delta_j$ such that
$\alpha_j'(x)\propto\alpha_j'(y)$ in $M$. As before, if
$x_{j+1}\in\alpha_{j+1}'(\Delta_{j+1}')$, we set
$M_{j+1}=\alpha_{j+1}'(\Delta_{j+1}')$,
$\alpha_{j+1}=\alpha_{j+1}'$, $\beta_j=\beta_j'$. Otherwise, let
$\Delta_{j+1}=\Delta_{j+1}'\oplus\langle e'\rangle$, where $e'$ is
an additional basis element. Define
${\alpha_{j+1}}_{|\Delta_{j+1}'}=\alpha_{j+1}'$,
$\alpha_{j+1}(e')=x_{j+1}$ and $\beta_j$ as the composition of
$\beta_j'$ with the natural inclusion of $\Delta_{j+1}'$ in
$\Delta_{j+1}$.

The proof follows then by induction.
\end{proof}

\noindent We are now ready to prove the main result of this
section.
\begin{thm}
\label{mainmonoid} Let $M$ be a countable refinement monoid. Then
there is a dimension monoid $\Delta$ and a surjective morphism
$\alpha\colon\Delta\to M$ such that $x\propto y$ for every pair of
elements $x$, $y$ in $\Delta$ that satisfy
$\alpha(x)\propto\alpha(y)$.
\end{thm}
\begin{proof}
Write $M=\{x_0,x_1,\ldots\}$. We take a commutative diagram of
simplicial monoids $\Delta_j$ and finitely generated submonoids
$M_j$ of $M$ as diagram~\ref{diagram} in
Proposition~\ref{almostthere}.

Define $\Delta=\varinjlim(\Delta_j,\beta_j)$, which is a dimension
monoid. The diagram (\ref{diagram}) induces a map
$\alpha\colon\Delta\to M$ by the universal property of inductive
limits. Denote by $\beta_{\infty,j}\colon\Delta_j\to\Delta$ the
inductive limit maps, and observe that we have a commutative
diagram
\[\begin{CD} \Delta_j @>{\alpha_j}>> M_j\\
@V{\beta_{\infty,j}}VV  @VV{\iota}V\\
\Delta @>{\alpha}>> M
\end{CD}\]
Note that
\[
\mathrm{Im}(\alpha)=\bigcup_j\alpha_j(\Delta_j)=\bigcup_jM_j=M\,,
\]
and so $\alpha$ is surjective.

Next, let $x$, $y\in\Delta$, and assume that
$\alpha(x)\propto\alpha(y)$. There exist then $j\geq 0$ and elements
$x_j$, $y_j$ in $\Delta_j$ such that $\beta_{\infty, j}(x_j)=x$
and $\beta_{\infty,j}(y_j)=y$. Now,
$\alpha_j(x_j)=\alpha(x)\propto\alpha(y)=\alpha_j(y_j)$, and hence
$\beta_j(x_j)\propto\beta_j(y_j)$ in $\Delta_{j+1}$, by
Proposition~\ref{almostthere}. But then
\[
x=\beta_{\infty,j}(x_j)=\beta_{\infty,j+1}(\beta_j(x_j))\propto\beta_{\infty,j+1}(\beta_j(y_j))=\beta_{\infty,j}(y_j)=y\,,
\]
and this finishes the proof.
\end{proof}

\begin{cor}[cf.\ \cite{gw}]
\label{cor:idealmonoid} Let $M$ be a countable
refinement monoid. There exists then a dimension monoid $\Delta$
such that $L(\Delta)\cong L(M)$.
\end{cor}

\begin{proof}
Follows immediately from Theorem~\ref{mainmonoid} and
Lemma~\ref{idealmonoid}.
\end{proof}

\noindent The preceding corollary was already known (although
maybe not with this precise formulation) and is related to a
result of G. M. Bergman about realizing distributive algebraic
lattices as ideal lattices of von Neumann regular rings
(\cite{berg}, see
also~\cite{werpubmat},~\cite{gw},~\cite{wtrans}). We make the
connection more explicit as follows (see also the comments after
Theorem~\ref{brat-ell}).

Let $M$ be a monoid. Define a congruence on $M$ by writing
$x\asymp y$ if $x\propto y$ and $y\propto x$. Set $\nabla (M)=M/\asymp\,$,
which in the literature is referred to as the \emph{maximal
semilattice quotient} of $M$ (see, e.g. \cite{cp}).  If $M$ is a
refinement monoid, then $\nabla (M)$ is also a refinement monoid,
see \cite[Lemma 2.4]{gw}. This is an example of a distributive
$0$-semilattice (which is, by definition, a commutative refinement
monoid such that $x+x=x$ for all $x$).

It was asked in \cite[Problem 10.1]{gw} (see also \cite[Problem
3]{tw}) whether for a distributive $0$-semilattice $S$ there is a
dimension group $G$ such that $\nabla (G^+)\cong S$ (where $G^+$
is the positive cone of $G$). The cases where $S$ is in fact a
lattice or where $S$ is countable were handled in \cite{gw}, with
positive solutions.

Let us note first that the countable case (\cite[Theorem 5.2]{gw})
can also be obtained from our results above. Namely, given a
countable distributive $0$-semilattice $S$, observe first that
$\nabla(S)\cong S$. Since $S$ is a refinement monoid, there exists
by Theorem~\ref{mainmonoid} a dimension monoid $N$ and a
surjective monoid morphism $\alpha\colon N\to S$ such that $x\propto
y$ whenever $x$ and $y$ in $N$ satisfy $\alpha(x)\propto\alpha(y)$.
This morphism clearly induces an isomorphism
$\nabla(N)\cong\nabla(S)$, given by $[x]\mapsto [\alpha(x)]$.
(Here $[x]$ denotes the class of the element $x$ modulo the
congruence $\asymp$\,.)

We now indicate how to prove Corollary~\ref{cor:idealmonoid} from
the methods in~\cite{gw}. We use first that if $M$ is a countable
refinement monoid, then $\nabla(M)$ is a countable distributive
$0$-semilattice. By~\cite[Theorem 5.2]{gw}, there is a dimension
monoid $N$ such that $\nabla(N)\cong\nabla(M)$, and
applying~\cite[Proposition 2.6 (iii)]{gw}, we find that $L(N)\cong
L(\nabla(N))\cong L(\nabla(M))\cong L(M)$, as desired.

In light of Theorem~\ref{mainmonoid}, it is also natural to ask
whether a similar result is available if the refinement monoid $M$
is not countable. The answer to this question is negative in
general. Specifically, there are semilattices $S$ of size
$\aleph_2$ that cannot be isomorphic to $\nabla(N)$ for any
dimension monoid $N$, as proved by R\r{u}\v zi\v cka (\cite{ru}).
Counterexamples for semilattices of size $\aleph_1$ were later
obtained by Wehrung in \cite{wtrans}.

\section{Applications to operator algebras}

\noindent
Let $A$ be a $C^*$-algebra. Recall (from Section~2) that $V(A)$ is the
monoid of Murray-von Neuman equivalence classes of projections in
$M_\infty(A)$, and that $V(A)$ is a refinement monoid if $A$
is of real rank zero. Define the \emph{dimension range} of $A$ to be
the subset $D(A)$ of $V(A)$ given by
\[
D(A)=\{[p]\in V(A)\mid p\mbox{ is a projection in }A\}\,.
\]
The subset $D(A)$ is a \emph{partial monoid} in the sense that
addition in $D(A)$ is only partially defined. Any $^*$-homomorphism
$\varphi \colon A \to B$ induces a monoid
morphism $V(\varphi) \colon V(A) \to V(B)$ defined by $[p] \mapsto
[\varphi(p)]$, and it satisfies $V(\varphi)(D(A))
\subseteq D(B)$.
If $A$ is unital, then an element $x$ in $V(A)$ belongs to $D(A)$ if
and only if $x \le [1_A]$. If $A$ is stable, then $D(A)=V(A)$.

The following two lemmas are well-known in the case where the target
algebra $A$ has the cancellation property (in which case these statements
are true with $V(A)$ replaced with $K_0(A)$; see for example
\cite[Lemma 7.3.2]{rll}).

\begin{lem}
\label{inducing0} Let $B$ be a finite dimensional $C^*$-algebra
and let $A$ be a unital $C^*$-algebra.
\begin{enumerate}
\item Suppose that $\varphi$, $\psi \colon B \to A$ are unital
$^*$-homomorphisms that satisfy $V(\varphi)=V(\psi)$. Then there
is a unitary $u$ in $A$ such that $u\varphi(x)u^*=\psi(x)$ for all
$x$ in $B$. \item Suppose that $\alpha \colon V(B) \to V(A)$ is a
monoid morphism
  that satisfies $\alpha([1_B]) = [1_A]$. Then there is a unital
  $^*$-homomorphism $\varphi \colon B \to A$ for which
  $V(\varphi) = \alpha$. If $e \in B$ and $f \in A$ are projections
  such that $\alpha([e]) = [f]$ and $\alpha([1_B-e]) = [1_A-f]$, then
  $\varphi$ above can be chosen satisfying $\varphi(e)
  = f$.
\end{enumerate}
\end{lem}

\begin{proof}
Write $B = M_{n_1} \oplus M_{n_2} \oplus \cdots \oplus M_{n_r}$ and let
$\{e_{ij}^{(k)}\}$ be a system of matrix units for $B$ (where
$1\leq k\leq r$ and $1\leq i,j\leq n_k$).

(i). Since
$\varphi(e_{11}^{(k)})\sim\psi(e_{11}^{(k)})$ for all $k$, we can
find partial isometries $v_k$ in $A$ such that
$v_k^*v_k=\varphi(e_{11}^{(k)})$ and
$v_kv_k^*=\psi(e_{11}^{(k)})$.
Set
$$u=\sum_{k=1}^r\sum_{i=1}^{n_k}\psi(e_{i1}^{(k)})v_k\varphi(e_{1i}^{(k)}).$$
One can now check that $u$ is a unitary element in $A$ with the
required
properties.

(ii). The matrix units $\{e_{ij}^{(k)}\}$ for $B$ can be chosen such
that $e = \sum_{(i,k) \in \Gamma} e_{ii}^{(k)}$ for some subset
$\Gamma$ of
$$\Omega := \{(i,k) \mid k = 1,2, \dots, r, \; i = 1,2, \dots, n_k\}.$$
We then have $1_B-e = \sum_{(i,k) \in \Gamma'} e_{ii}^{(k)}$, when
$\Gamma' = \Omega \setminus \Gamma$.

Find, for some
large enough $m$, pairwise
orthogonal projections $g_i^{(k)}$ in $M_m(A)$ such
that $[g_i^{(k)}] = \alpha([e_{ii}^{(k)}])$ for all $(i,k) \in
\Omega$. Put $g =
\sum_{(i,k) \in \Gamma} g_i^{(k)}$ and $g' =
\sum_{(i,k) \in \Gamma'} g_i^{(k)}$. Then
$$
[g] = \sum_{(i,k) \in \Gamma} [g_i^{(k)}] = \sum_{(i,k) \in \Gamma}
\alpha([e_{ii}^{(k)}]) = \alpha\big(\big[\sum_{(i,k) \in \Gamma}
e_{ii}^{(k)}\big]\big) = \alpha([e]) = [f],$$
and, similarly, $[g'] = [1_A-f]$.
It follows that
$$g=v^*v, \quad f = vv^*, \quad g' = w^*w, \quad 1_A-f = ww^*,$$
for some partial isometries $v$ and $w$ in
$M_{1m}(A)$. Put $u = v+w$ and put $f_{ii}^{(k)} = ug_i^{(k)}u^*$ for
all $(i,k)$ in $\Omega$. Then $\{f_{ii}^{(k)}\}$ are
pairwise orthogonal projections in $A$,
$$[f_{ii}^{(k)}] = [g_{i}^{(k)}] = \alpha([e_{ii}^{(k)}]), \quad
\sum_{(i,k) \in \Gamma} f_{ii}^{(k)} = f, \quad \sum_{(i,k) \in \Gamma'}
f_{ii}^{(k)} = 1_A-f.$$

Now, $[f_{ii}^{(k)}]= \alpha([e_{ii}^{(k)}]) =\alpha([e_{jj}^{(k)}])
=[f_{jj}^{(k)}]$, that is, $f_{ii}^{(k)}\sim f_{jj}^{(k)}$ for fixed
$k$ and for
all $i$ and $j$. It is standard (see for example \cite[Lemma
7.1.2]{rll}) that the system $\{f_{ii}^{(k)}\}$ extends to a system of
matrix units $\{f_{ij}^{(k)}\}$ for $B = M_{n_1} \oplus M_{n_2} \oplus
\cdots \oplus M_{n_r}$, and that there is a $^*$-homomorphism $\varphi
\colon B \to A$ given by $\varphi(e_{ij}^{(k)}) = f_{ij}^{(k)}$.

We must check that $\varphi$ has the desired properties, and note first
$$\varphi(e) = \varphi\big(\sum_{(i,k) \in \Gamma} e_{ii}^{(k)} \big)
= \sum_{(i,k) \in \Gamma} f_{ii}^{(k)} = f,$$
and, similarly, $\varphi(1_B-e) = 1_A-f$. In particular, $\varphi(1_B)
= 1_A$. Finally, as $[e_{ii}^{(k)}]$ generate $V(B)$ and
$$V(\varphi)([e_{ii}^{(k)}]) = [\varphi([e_{ii}^{(k)}])] =
[f_{ij}^{(k)}] = \alpha([e_{ii}^{(k)}]),$$
we conclude that $V(\varphi) = \alpha$.
\end{proof}

\begin{lem}
\label{inducing-s} Let $A$ be a stable $C^*$-algebra, let $B$ be an
AF-algebra, and let $\alpha \colon V(B) \to V(A)$ be a monoid
morphism. Then there is a $^*$-homomorphism $\varphi \colon B\to A$ such that
$V(\varphi)=\alpha$.
\end{lem}

\begin{proof}
The proof consists of constructing the commutative diagram:
\begin{equation} \label{diagram2}
\begin{split}
\xymatrix{B_1 \ar[r]^{\psi_1} \ar[d]_{\varphi_1} & B_2
  \ar[r]^{\psi_2} \ar[d]_{\varphi_2} &  B_3 \ar[r]^{\psi_3}
  \ar[d]_{\varphi_3} &  \cdots \ar[r] & B \ar@{-->}[d]^\varphi \\
f_1Af_1 \ar[r]^\iota & f_2Af_2 \ar[r]^\iota & f_3Af_3
\ar[r]^\iota & \cdots \ar[r] & A_0}
\end{split}
\end{equation}
where $B$ is written as an inductive limit of finite dimensional
$C^*$-algebras $B_n$ with connecting maps $\psi_n$ (that are not
necessarily unital). The projections $f_1, f_2, \dots$, constructed below,
form an increasing sequence, and they satisfy $[f_n] = \alpha \big(
[\psi_{\infty,n}(1_{B_n})]\big)$, where $\psi_{\infty,n} \colon B_n \to B$
is the inductive limit map. 
The algebra $A_0$ in \eqref{diagram2} is the hereditary sub-\Cs{} of
$A$ given by
$$A_0 = \overline{\bigcup_{n=1}^\infty f_nAf_n}.$$

We proceed to find the projections $f_n$. Put $e_n =
\psi_{\infty,n}(1_{B_n})$. Then $e_1,e_2, \dots$ is an increasing
sequence of projections in $B$. Since
$A$ is stable we can find pairwise orthogonal projections
$g_1,g_2, \dots$ in $A$ such that $[g_1] = \alpha([e_1])$ and $[g_n] =
\alpha([e_n - e_{n-1}])$ for $n \ge 2$.
Put $f_n = g_1+g_2+ \cdots+g_n$. Then $f_1,f_2, \dots$ is an
increasing sequence of projections in $A$ and
$[f_n] = \alpha \big([\psi_{\infty,n}(1_{B_n})]\big) = (\alpha \circ
V(\psi_{\infty,n}))([1_{B_n}])$.

The next step is to find unital $^*$-homomorphisms $\varphi_n \colon
B_n \to f_nAf_n$ making 
Diagram~\eqref{diagram2} commutative and such that
$V(\varphi_n) = \alpha \circ V(\psi_{\infty,n})$. The existence
of $\varphi_1$ follows from Lemma~\ref{inducing0}~(ii) (with
$e=f=0$). At step two, use again Lemma~\ref{inducing0}~(ii) to find a
unital $^*$-homomorphism $\psi_2 \colon B_2 \to f_2Af_2$ such that
$V(\psi_2) = \alpha \circ V(\psi_{\infty,2})$ and
$\psi_2(\varphi_1(1_{B_1})) = f_1$. The two $^*$-homomorphisms
$\psi_2 \circ \varphi_1$ and $\iota \circ \varphi_1$ co-restrict to
unital $^*$-homomorphisms $B_1 \to f_1Af_1$, and $V(\psi_2 \circ
\varphi_1) = V(\iota \circ \varphi_1)$. Use Lemma~\ref{inducing0}~(i) to
find a unitary $u_0$ in $f_1Af_1$ such that ${\mathrm{Ad}}_{u_0} \circ \psi_2
\circ \varphi_1 = \iota \circ \varphi_1$. Put $u = u_0 +
(f_2-f_1)$ (which is a unitary in $f_2Af_2$) and put $\varphi_2 =
{\mathrm{Ad}}_u \circ \psi_2$. Then $\varphi_2$ is a unital $^*$-homomorphism,
$V(\varphi_2) = \alpha \circ V(\psi_{\infty,2})$, and the first
square in the diagram \eqref{diagram2} commutes.

Continuing in this way we find the remaining $^*$-homomorphisms
$\varphi_n$ making the successive squares in diagram \eqref{diagram2}
commutative. By the 
universal property of inductive limits there is a
$^*$-homomorphism $\varphi \colon B \to A$ that makes \eqref{diagram2}
commutative. For $x$ in $V(B_n)$ one has
$$V(\varphi)\big(V(\psi_{\infty,n})(x)\big) = V(\varphi \circ
\psi_{\infty,n})(x) =
V(\varphi_n)(x) = \alpha \big(V(\psi_{\infty,n})(x)\big).$$
The functor $V$ is continuous, which entails that $V(B) =
\cup_{n=1}^\infty V(\psi_{\infty,n})(V(B_n))$. We can now conclude that
$V(\varphi) = \alpha$.
\end{proof}

\noindent At first sight it is tempting to believe that Lemma~\ref{inducing-s}
holds without the assumption that $A$ is stable, but 
with the extra
assumption $\alpha(D(B)) \subseteq D(A)$.  
At second thought, it appears 
plausible (to the authors) that there might exist a
(non-stable)
$C^*$-algebra $A$ for which there is a monoid morphism $\alpha \colon
V({\mathcal{K}}) \to V(A)$, that satisfies $\alpha(D({\mathcal{K}})) \subseteq
D(A)$, and that does not lift to a $^*$-homomorphism $\varphi \colon
{\mathcal{K}} \to A$.


\begin{thm}
\label{embedding} Let $A$ be a separable $C^*$-algebra of real
rank zero. Then there exists an $AF$-algebra $B$ and a
$^*$-monomorphism $\varphi \colon B \to A$ such that
\begin{enumerate}
\item $\varphi$ induces an isomorphism of ideal lattices,
\item for any two projections $e,f$ in $B$ one has $[e] \propto [f]$ in
  $V(B)$ if and only if $[\varphi(e)] \propto [\varphi(f)]$ in $V(A)$,
\item $V(\varphi) \colon V(B) \to V(A)$ is onto. 
\item $A = \overline{\varphi(B)A\varphi(B)}$.
\end{enumerate}
If $A$ is unital, then the AF-algebra $B$ is necessarily also unital,
and $\varphi$ is unit preserving.
\end{thm}

\begin{proof}
Since $A$ has real rank zero, we know that $V(A)$ is a refinement
monoid (\cite[Lemma~2.3]{ap}). It is also countable because $A$ is
separable. Therefore we can apply Theorem~\ref{mainmonoid} to find
a dimension monoid $\Delta$ and a surjective monoid morphism
$\alpha \colon \Delta \to V(A)$ such that $x\propto y$ for
every pair of elements $x$, $y$ in $\Delta$ for which
$\alpha(x) \propto \alpha(y)$.

By the structure theory for AF-algebras (see e.g.\ \cite{eff}),
there is a stable AF-algebra $B_s$ such that $V(B_s)$ is isomorphic to
$\Delta$. Identifying $\Delta$ with $V(B_s)$ we will assume that the
domain of $\alpha$ is $V(B_s)$. Let
$A_s$ denote the stabilization of $A$ and identify $A$ with a (full) hereditary
sub-$C^*$-algebra of $A_s$. Identify $V(A)$ and $V(A_s)$, so that
$\alpha$ is identified with a monoid morphism $V(B_s) \to V(A_s)$. Use
Lemma~\ref{inducing-s} to lift $\alpha$ to a $^*$-homomorphism $\psi \colon B_s
\to A_s$.

Choose an increasing approximate unit $\{f_n\}_{n=1}^\infty$
consisting of projections for $A$ (so that $A =
\overline{\cup_{n=1}^\infty f_n A_s f_n}$). By surjectivity of
$V(\psi)$, and since $B_s$ is stable, there are pairwise orthogonal
projections $g_1,g_2, \dots$ in $B_s$ such that
$[\psi(g_1)] = [f_1]$ and $[\psi(g_n)] = [f_n-f_{n-1}]$ for $n \ge
2$. Put $e_n = g_1+g_2+ \cdots + g_n$. Then $e_1,e_2, \dots$ is an
increasing sequence of projections in $B_s$ and $[\psi(e_n)] =
[f_n]$. Put $B = \overline{\cup_{n=1}^\infty e_n B_s e_n}$, so that
$B$ is a hereditary sub-$C^*$-algebra of $B_s$, and hence is an
AF-algebra. We claim that $B$ is
full in $B_s$.
Indeed, if $p$ is a projection in $B_s$, then $\alpha([p]) \propto
[f_n] = \alpha([e_n])$ for some $n$, because $A$ is full in $A_s$, and
this implies  $[p] \propto [e_n]$ by the special property of
$\alpha$. Thus $p$ belongs to the closed two-sided ideal in $B_s$
generated by $B$; and as $p$ was arbitrary, $B$ is full in $B_s$.

Choose partial isometries $v_n$ in $A$ such that
$$v_1v_1^* = \psi(g_1), \quad v_1^*v_1 = f_1, \qquad v_nv_n^* =
\psi(g_n), \quad v_n^*v_n = f_n - f_{n-1}, \; n \ge 2,$$
set $u_n = v_1+v_2+ \cdots + v_n$, so that $u_nu_n^* = \psi(e_n)$ and
$u_n^*u_n = f_n$, and define $\varphi_n \colon
e_nB_se_n \to f_nA_sf_n$ by $\varphi_n(b) = u_n^*\psi(b)u_n$. Then
  we obtain a commutative diagram
\begin{equation} \label{diagram3}
\begin{split}
\xymatrix{e_1B_se_1 \ar[d]^{\varphi_1} \ar[r]^\iota & e_2B_se_2
  \ar[d]^{\varphi_2} \ar[r]^\iota & e_3B_se_3 \ar[d]^{\varphi_3}
  \ar[r]^\iota & \cdots \ar[r] & B  \ar@{-->}[d]^\varphi \\
f_1A_sf_1 \ar[r]^\iota & f_2A_sf_2 \ar[r]^\iota & f_3A_sf_3 \ar[r]^\iota &
\cdots \ar[r] & A}
\end{split}
\end{equation}
which by the universal property of inductive limits produces the
$^*$-homomorphism $\varphi$. Since $\varphi(e_n) = \varphi_n(e_n) =
f_n$ we see that (iv) holds.

We claim that $V(\varphi) = \alpha$. Indeed,
$V(\varphi_n) = V(\psi|_{e_nB_se_n})$. Hence, if we apply the continuous
functor $V$ to the diagram \eqref{diagram3}, then the
resulting diagram remains commutative if $V(\varphi)$ is replaced with
$V(\psi|_B)$, whence $V(\varphi) = V(\psi|_B)$ by the universal
property of inductive limits. Identifying $V(A)$ and $V(B)$ with
$V(A_s)$
and $V(B_s)$, respectively, we obtain $V(\varphi) = V(\psi) =
\alpha$.

Now, (ii) follows by the stipulated property of $\alpha$, and (i)
follows from Corollary~\ref{ideal2}, part (ii), and the assumption that
$A$ is of real rank zero. It follows in particular from (i) that
$\varphi$ is injective.

Surjectivity of $V(\varphi) \colon V(B) \to V(A)$ follows from the
surjectivity of $\alpha$.

Suppose finally that $A$ is unital. Then $1_A$ belongs to $\varphi(B)$
by (iv). This entails that $\varphi(B)$ is unital with unit
$1_A$. Because $\varphi$ is injective, and hence an isomorphism from
$B$ to $\varphi(B)$, it follows that $B$ is unital and that $\varphi$
maps the unit of $B$ onto the unit of $\varphi(B)$, i.e.
$\varphi(1_B)=1_A$, as desired.
\end{proof}

\begin{cor} \label{cor:embedding}
Let $A$ be a simple, separable $C^*$-algebra of real
rank zero. Then there is a simple AF-algebra $B$ and an
embedding $\varphi \colon B \to A$ such that $V(\varphi) \colon V(B)
\to V(A)$ is surjective, and $A
= \overline{\varphi(B)A\varphi(B)}$.
\end{cor}

\noindent One of our motivations for pursuing Theorem~\ref{embedding}
was to establish---what later turned out to be a well-known fact!
(see Theorem~\ref{brat-ell} below)---that
the primitive ideal space of an arbitrary, 
separable, real rank zero algebra is
homeomorphic to the primitive ideal space of some
AF-algebra. Combining this with Kirchberg's seminal
classification of $\mathcal{O}_2$-absorbing separable, nuclear
$C^*$-algebras (see Theorem~\ref{Kirchberg} below) one obtains
Corollary~\ref{cor:O_2}, that the tensor product of an arbitrary
separable, stable, nuclear $C^*$-algebra of real rank zero with
$\mathcal{O}_2$ is isomorphic to an AF-algebra tensor $\mathcal{O}_2$.

Recall that the primitive ideal space, $\mathrm{Prim}(A)$, of a
$C^*$-algebra $A$ is the set of all kernels of irreducible
representations of $A$ equipped with the hull-kernel topology
(also known as the Jacobson topology). There is a bijective
correspondence
\[
I\mapsto\mathrm{hull}(I)=\{P\in\mathrm{Prim}(A)\mid I\subseteq P\}
\]
between the set of closed ideals of $A$ and the set
of closed subsets of $\mathrm{Prim}(A)$ (so $\mathrm{Prim}(A)$
contains the same information as $L(A)$, the ideal lattice of $A$). It
is known that $\mathrm{Prim}(A)$ is a locally compact
$T_0$-space with the Baire property, but there is no description of
which spaces with these properties arise as $\mathrm{Prim}(A)$ for
some $C^*$-algebra $A$.

In the real rank zero
case, or more generally, in the case of  $C^*$-algebras with property (IP)
(each ideal in the $C^*$-algebra is generated by its projections) we
have the following theorem of Bratteli and Elliott, \cite{brael}.
We write $X \cong Y$ when $X$ and $Y$ are homeomorphic topological spaces.

\begin{thm}[Bratteli--Elliott] \label{brat-ell}
Let $X$ be a locally compact $T_0$-space with the Baire property and
with a countable basis. Then
the following conditions are equivalent.
\begin{enumerate}
\item $X$ has a basis of compact-open sets.
\item $X \cong {\mathrm{Prim}}(A)$ for some AF-algebra $A$.
\item $X \cong {\mathrm{Prim}}(A)$ for some $C^*$-algebra $A$ with
  property (IP).
\end{enumerate}
\end{thm}

\noindent It follows either from Theorem~\ref{brat-ell} or from
our Theorem~\ref{embedding} that if $A$ is any separable
$C^*$-algebra of real rank zero, then ${\mathrm{Prim}}(A) \cong
{\mathrm{Prim}}(B)$ for some AF-algebra $B$. Alternatively, this fact 
follows from our Corollary~\ref{cor:idealmonoid} combined Zhang's
theorem that the ideal lattice of a $C^*$-algebra $A$ of real rank
zero is (canonically) isomorphic to the ideal lattice of $V(A)$.

The theorem below was proved by Kirchberg, \cite{Kir:class}.

\begin{thm}[Kirchberg] \label{Kirchberg}
Let $A$ and $B$ be separable, nuclear $C^*$-algebras. Then
$${\mathrm{Prim}}(A) \cong {\mathrm{Prim}}(B) \iff A \otimes
{\mathcal{O}}_2 \otimes {\mathcal{K}} \cong B \otimes
{\mathcal{O}}_2 \otimes {\mathcal{K}}\,.$$
\end{thm}

\begin{lem} \label{lm:O_2-tensor-B}
Let $B$ be a separable $C^*$-algebra of real rank zero. Then every
projection in $B 
\otimes \mathcal{O}_2$ is equivalent to a projection in $B \otimes \C 1
\subseteq B \otimes \mathcal{O}_2$.
\end{lem}

\begin{proof} Let $p$ be a projection in $B \otimes \mathcal{O}_2$, and let
  $J$ be the closed two-sided ideal in $B \otimes \mathcal{O}_2$ generated by
  $p$. As $\mathcal{O}_2$ is nuclear it follows from 
\cite[Theorem 3.3]{Bla:tensor}
  that $J = J_0
  \otimes \mathcal{O}_2$ for some closed two-sided ideal $J_0$ in $B$. We
  proceed to show that $J_0$ contains a full projection. As $B$ is
  separable and of real rank zero, $J_0$ has an increasing approximate
  unit $\{e_n\}_{n=1}^\infty$ consisting of projections. Let $I_n$ be
  the closed two-sided ideal in $J_0$ generated by $e_n$. Then
$$p \in J=\overline{\bigcup_{n=1}^\infty I_n \otimes \mathcal{O}_2}\,,$$
whence $p$ belongs to $I_n \otimes \mathcal{O}_2$ for some $n$ (because $p$ is
a projection). Put $q = e_n \otimes 1 \in (B \otimes \C 1) \cap
J$. Then $p$ and $q$ are full projections in $J$; and they are
properly infinite because $J$ is purely infinite (cf.\
\cite[Proposition~4.5 and Theorem~4.16]{KirRor:pi}). As $K_0(J) = 0$
(see eg.\ \cite[Theorem~2.3]{Cuntz:KOn}) it 
follows from \cite[Section~2]{Cuntz:KOn} that $p$ and $q$ are
equivalent. 
\end{proof}

\noindent Combining Kirchberg's theorem with the result about the
primitive ideal space of real rank zero algebras mentioned above yields:

\begin{cor} \label{cor:O_2}
 Let $A$ be a separable, nuclear $C^*$-algebra of real rank
  zero. Then $A \otimes {\mathcal{O}}_2$ is of real rank zero, 
  there is an AF-algebra 
  $B$, and there is a sequence of natural numbers
  $\{r_n\}_{n=1}^\infty$ such that
$$A \otimes {\mathcal{O}}_2  \, \cong \, B  \otimes
{\mathcal{O}}_2 \; \cong \; \lim_{n \to \infty}
\bigoplus_{j=1}^{r_n}  {\mathcal{O}}_2\,.$$
\end{cor}

\begin{proof} The right-most $C^*$-algebra 
displayed above is of real rank zero, so
  the first claim follows from the last claim. If $B$ is an
  AF-algebra, then $B = 
  \lim_{n \to \infty} B_n$, where $B_n$ is the direct sum of $r_n$ 
  full matrix algebras. As $M_k \otimes {\mathcal{O}}_2$ is isomorphic to
  ${\mathcal{O}}_2$ for all $k$, we see that $B_n \otimes
  {\mathcal{O}}_2$ is isomorphic to $\oplus_{j=1}^{r_n}
  {\mathcal{O}}_2$. This proves the second isomorphism.

Use Theorem~\ref{brat-ell} or Theorem~\ref{embedding} to
  find a stable AF-algebra $B_s$ with ${\mathrm{Prim}}(A) \cong
  {\mathrm{Prim}}(B_s)$. Then $A \otimes {\mathcal{O}}_2 \otimes
  {\mathcal{K}}$ is isomorphic to $B_s  \otimes
{\mathcal{O}}_2$ by Kirchberg's Theorem~\ref{Kirchberg}. This implies
that $A \otimes {\mathcal{O}}_2$ is isomorphic to a hereditary
sub-$C^*$-algebra $D$ of $B_s  \otimes {\mathcal{O}}_2$. As $B_s  \otimes
{\mathcal{O}}_2$ is separable and of real rank zero, there is an
increasing approximate unit $\{q_n\}_{n=1}^\infty$ consisting of
projections for $D$. Use Lemma~\ref{lm:O_2-tensor-B} and stability of
$B_s$ to find pairwise orthogonal 
projections 
$f_1,f_2,\dots $ in $B_s$ such that
$f_1 \otimes 1 \sim q_1$ and $f_n \otimes 1 \sim q_n - q_{n-1}$ for $n
\ge 2$. Choose partial isometries 
$v_n$ in $B_s  \otimes
{\mathcal{O}}_2$ such that $v_n^*v_n = f_n \otimes 1$, $v_nv_n^* =
q_n-q_{n-1}$, set 
$u_n = v_1+v_2 + \cdots +v_n$, set $p_n = f_1+f_2 + \cdots +
f_n$ and set $B = \overline{\cup_{n=1}^\infty p_n B_s p_n}$. Then $B$
is an AF-algebra, $u_n^*u_n = p_n \otimes 1$, $u_nu_n^* = q_n$, and we
have a commutative diagram 
$$
\xymatrix{q_1( {\mathcal{O}}_2 \otimes B_s)q_1 \ar[d]^{\varphi_1}
  \ar[r]^\iota & q_2( {\mathcal{O}}_2 \otimes B_s)q_2
  \ar[d]^{\varphi_2} \ar[r]^\iota &  q_3( {\mathcal{O}}_2 \otimes B_s)q_3
 \ar[d]^{\varphi_3}
  \ar[r]^\iota & \cdots \ar[r] & D  \ar@{-->}[d]^\varphi \\
(p_1B_sp_1) \otimes  {\mathcal{O}}_2 
 \ar[r]^\iota & (p_2B_sp_2) \otimes  {\mathcal{O}}_2 
 \ar[r]^\iota & (p_3B_sp_3) \otimes  {\mathcal{O}}_2 
 \ar[r]^-\iota &
\cdots \ar[r] & B \otimes  {\mathcal{O}}_2}
$$
where $\varphi_n(x)= u_n^*xu_n$. The  $^*$-homomorphism $\varphi$
induced by the diagram is an isomorphism because each $\varphi_n$ is
an isomorphism. It follows that $A  \otimes {\mathcal{O}}_2 \cong D
\cong  B \otimes  {\mathcal{O}}_2$.
\end{proof}

\noindent In \cite{lin2}, Lin considered a class $\mathcal{A}$ of
separable $C^*$-algebras that have 
trivial $K$-Theory. This class consists of those $C^*$-algebras $A$
that can be written as inductive 
limits of finite direct sums of hereditary sub-$C^*$-algebras of
$\mathcal{O}_2$ (\cite[Corollary 
3.11]{lin2}), and can be completely classified by their monoids of
equivalence classes of projections 
(\cite[Theorem 3.13]{lin2}). Since every hereditary sub-$C^*$-algebra
of $\mathcal{O}_2$ is isomorphic to 
$\mathcal{O}_2$ (in the unital case), or to
$\mathcal{O}_2\otimes\mathcal{K}$ (in the non-unital case),  
it follows from Corollary~\ref{cor:O_2}, that the $C^*$-algebras in
the class $\mathcal{A}$ are precisely those of the form
$A\otimes\mathcal{O}_2$, where $A$ is a (separable) AF-algebra. In the 
unital case, this was already observed in \cite[Theorem 9.4]{gw}. The
methods used there were based on  
the fact that for any $A\in\mathcal{A}$, the monoid $V(A)$ is a
distributive $0$-semilattice (cf. \cite{lin2}). 

The range of the invariant for the $C^*$-algebras in the class
$\mathcal{A}$ was also described in 
\cite[Section 4]{lin2}, as follows. Given any countable distributive
$0$-semilattice $V$ with a largest 
element, there is a unital $C^*$-algebra $A$ in $\mathcal{A}$ such
that $V(A)\cong V$ (\cite[Theorem 
4.14]{lin2}). We note that our methods developed in Section~3,
combined with the observations above, 
allow us to recover this result. Namely, given a countable
distributive $0$-semilattice $V$, find by 
Theorem~\ref{mainmonoid} a dimension monoid $N$ and a surjective
monoid morphism $\alpha\colon N\to V$ 
such that $\alpha(x)\propto\alpha(y)$, whenever $x\propto y$. As in
the comments following Theorem ~\ref{mainmonoid}, this induces an
isomorphism $\nabla(N)\cong\nabla(V)\cong V$ (where $\nabla(\cdot)$
denotes the passage to the maximal semilattice quotient). By the
structure theory of AF-algebras, there is a stable AF-algebra $B$ such
that $V(B)\cong N$. But now $V(B\otimes\mathcal{O}_2) \; \cong \;
V(B)/\!\asymp \; = \; \nabla(V(B)) \; \cong \; \nabla(N)$, and from
this it follows that $V(B\otimes\mathcal{O}_2) \, \cong \, \nabla(V)
\, \cong  \, V$. 

In the case that $V$ has a maximal element $u$, this is an order-unit
of the monoid $V$, and the surjectivity of the morphism $\alpha$ in
the paragraph above provides an element $v$ in $N$ such that
$\alpha(v)=u$ and that will be an order-unit for $N$. In this case,
the AF-algebra $B$ can be chosen to be unital.

One can combine Bratteli and Elliott's Theorem~\ref{brat-ell}
with Kirchberg's Theorem~\ref{Kirchberg} to conclude that $A \otimes
{\mathcal{O}}_2$ is stably isomorphic
to an AF-algebra tensor $\mathcal{O}_2$ whenever $A$ is a separable
nuclear $C^*$-algebra with property (IP). In particular, for such
$C^*$-algebras $A$, the tensor product $A \otimes
{\mathcal{O}}_2$ is of real rank zero; a curious fact that we
expect is derivable by more direct means. (Warning: There are
separable nuclear $C^*$-algebras $A$ with the ideal property where $A
\otimes \mathcal{O}_\infty$ is not of real rank zero!)

\section{Divisible $C^*$-algebras}

\noindent Divisibility in the context of a monoid $M$ refers to the
property that the equation $nx = y$ has a solution $x \in M$ 
(in which
case we say that $n$ divides $y$ in $M$) for all
$y$ in $M$ and for all (or some) large natural numbers $n$. If $A$ is a
$C^*$-algebra and $p$ is a projection in $A$, then $nx = [p]$ has a
solution $x \in V(A)$ if and only if there is a unital embedding $M_n
\to pAp$. Divisibility in this strong form is rare. 
A weaker
form of divisibility, which is much more frequent --- and
still useful --- is described in the definition below.

\begin{definition} \label{def:1}
Let $A$ be a $C^*$-algebra and let $p$ be a non-zero
  projection in $A$. We say that $A$ is \emph{weakly divisible of
    degree $n$ at $p$} if there is a unital $^*$-homomorphism
\begin{equation} \label{eq:embedding}
M_{n_1} \oplus M_{n_2} \oplus \cdots \oplus M_{n_r} \to pAp
\end{equation}
for \emph{some} natural numbers $r, n_1, n_2, \dots, n_r$ where $n_j \ge n$
for all $j$. If there is a unital $^*$-homomorphism as in
\eqref{eq:embedding} for \emph{each} set $r, n_1,n_2, \dots, n_r$
of natural numbers for which
  $\gcd(n_1,n_2, \dots, n_r)$ divides $[p]$ in $V(A)$, then we say
  that $A$ is \emph{weakly divisible at $p$}. Finally, if $A$ is
  weakly divisible at $p$ for every non-zero projection $p$ in $A$,
  then $A$ is called \emph{weakly divisible}.
\end{definition}

\noindent 
The notion of \emph{approximate
  divisibility} was introduced in \cite{BlaKumRor}. We recall the
definition: A unital 
$C^*$-algebra $A$ is approximately divisible if there is a sequence of
unital $^*$-homomorphisms $\varphi_n \colon M_2 \oplus M_3 \to A$ such
that $\varphi_n(x)a-a\varphi_n(x) \to 0$ for all $x \in M_2 \oplus
M_3$ and all $a \in A$. Being approximately divisible implies being
weakly divisible (at all projections $p$ in $A$) (see
\cite{BlaKumRor}). The crucial difference between weak divisibility
and approximate divisibility  is
the assumption of asymptotic centrality in the latter. Approximately
divisible $C^*$-algebras are very well behaved. In
particular, any simple, approximately divisible $C^*$-algebra is
either stably finite or purely infinite, and its ordered $K_0$-group
is always weakly unperforated (see \cite{BlaKumRor}).

We show here that weak divisibility is almost automatic for
$C^*$-algebras of real rank zero. There are examples of
non-nuclear, simple, unital $C^*$-algebras of
real rank zero that are 
weakly divisible but not approximately divisible
(\cite{DykRor}). Perhaps the most fundamental---and optimistic!---open
question concerning $C^*$-algebras of real rank zero is the following:

\begin{question} \label{apprdiv} Is every simple, unital, nuclear,
  non-type I $C^*$-algebra of real rank zero approximately divisible?
\end{question}

\noindent
If $A$ is weakly divisible at $p$, then there is a
  unital $^*$-homomorphism from $M_{n_1} \oplus M_{n_2} \oplus \cdots
  \oplus M_{n_r}$ into $pAp$ for every finite set
$n_1,n_2, \dots, n_r$ of natural numbers such that
  $\gcd(n_1,n_2, \dots, n_r)=1$.

There is a similar notion of weak divisibility in a monoid $M$:
$M$ is weakly 
divisible of degree $n$ at an element $x \in M$ if
there are natural numbers $r, n_1, \dots, n_r$ and elements $x_1,
\dots, x_r$ in $M$ such that $n_j \ge n$ for all $j$ and $x =
n_1x_1 + \cdots + n_rx_r$. It is easily seen that a $C^*$-algebra $A$
is weakly divisible of degree $n$ at a projection $p \in A$ if and
only if $V(A)$ is weakly divisible of degree $n$ at $[p]$.

If $A$ is weakly divisible, then so is any quotient, ideal, and
corner of $A$.

Since $M_2 \oplus M_3$ maps unitally (but not necessarily injectively)
into any matrix algebra $M_n$, with $n \ge 2$, weak divisibility of
degree 2 at a projection $p$ means that there is a unital
$^*$-homomorphism $M_2 \oplus M_3 \to pAp$.

Clearly, if $A$ is weakly divisible at a projection $p$,
then $A$ is
weakly divisible of degree 2 at $p$. In the converse direction, we have the
following:

\begin{lem} \label{divisible1}
Let $A$ be a $C^*$-algebra. If $A$ is weakly divisible of degree 2 at
\emph{all}
(non-zero) projections in $A$, then $A$ is weakly divisible.
\end{lem}

\begin{proof} We prove this in the monoid language. Hence we assume
  that for every $x$ in $V(A)$ there are elements $y,z$ with $x
  =2y+3z$.
Let now $x$ in $V(A)$ be fixed and let $n_1, \dots, n_r$ be
  natural numbers such that $d = \gcd(n_1, \dots, n_r)$ divides
  $x$. It suffices to consider the case where
$d=1$. Indeed, $x = dy$ for some $y$ in
$V(A)$. Put $m_j = n_j/d$. Then $\gcd(m_1, \dots, m_r)=1$, and if there
are elements $y_1, \dots, y_r$ in $V(A)$ with $y = m_1y_1 + \dots
m_ry_r$, then $x = n_1y_1+ \cdots +n_ry_r$.

Assume now that $d=1$. There is a natural number $m_0$ such that all
natural numbers $m \ge m_0$ belong to the sub-semigroup of the natural
numbers generated by $n_1, \dots, n_r$. Let $k$ be a natural number
such that $2^k \ge m_0$. By $k$ successive applications of
weak divisibility (first to $x = 2y+3z$, second to $y$ and $z$, and so
forth) we obtain elements $y_0,y_1, \dots, y_k$ in $V(A)$ such that
$$x = 2^k y_0 + 2^{k-1}3 y_1 + \cdots + 3^k y_k.$$
Since $2^{k-l} \! \cdot \! 3^{l} \ge m_0$ we can write $2^{k-l} \!
\cdot \! 3^{l} = \sum_{j=1}^r d_{l,j} \, n_j$
for suitable non-negative integers $d_{l,j}$. Now,
$$x = \sum_{l=0}^k (2^{k-l} \! \cdot \! 3^l) y_l =  \sum_{l=0}^k
\Big( \sum_{j=1}^r d_{l,j} \, n_j \Big) y_l = \sum_{j=1}^r n_j \Big(
\sum_{l=0}^r d_{l,j}\, y_l \Big) = \sum_{j=1}^r n_j x_j,$$
as desired, when $x_j = \sum_{l=0}^k d_{l,j}\, y_l$.
\end{proof}

\noindent Our main results on divisibility in \Cs s of real rank zero
are contained in Propositions~\ref{prop:div1} and
\ref{prop:divisible2} and in Theorem~\ref{thm:divisibility}
below. Proposition~\ref{prop:div1} is actually a special case of
Theorem~\ref{thm:divisibility}, but is emphasized because of its
independent interest and because it is used in the proof of the more general
results. 

A simple separable $C^*$-algebra $A$ is of type I precisely when
it is isomorphic to a matrix algebra $M_n$, for some $n$, or to
the compact operators, $\mathcal{K}$, on some separable Hilbert
space $H$. In other words, $A$ is of type I precisely when it is
isomorphic to a sub-$C^*$-algebra of $\mathcal{K}$.

\begin{prop} \label{prop:div1}
Every  simple, separable \Cs{} of real rank zero, that is not of type I,
is weakly divisible.  
\end{prop}

\begin{proof} By Lemma~\ref{divisible1} it suffices to show that $A$
  is weakly divisible of degree 2 at every non-zero projection $p$ in
  $A$. By Theorem~\ref{embedding} 
  (see also Corollary~\ref{cor:embedding}) there is
  a simple unital AF-algebra 
  $B$ and a unital embedding $\varphi \colon B \to pAp$ such that
  $V(\varphi) \colon V(B) \to V(pAp)$ is onto. If $B$ is weakly divisible of
  degree 2 at its unit, then $A$ is weakly divisible of degree 2 at $p$. 
  This is well-known to be the case
  when $B$ is infinite dimensional (see for example
  \cite[Proposition~4.1]{BlaKumRor}) and it is trivially true
  when $B \cong M_n(\C)$ for $n \ge 2$. 

Suppose now that $B= \C$ (the only case where $B$ is not weakly divisible of
degree 2 at its unit). We show that $pAp$ is properly infinite, and
this will imply that the Cuntz algebra $\mathcal{O}_\infty$, and hence $M_2
\oplus M_3$, embed unitally into $pAp$, thus showing that $A$ is
weakly divisible of degree 2 at $p$ also in this case.

The assumption that $A$ is not of type I (and $p \ne 0$) implies that
$pAp \ne \C p$. The corner $pAp$ therefore contains a non-trivial
projection $q$. By surjectivity of $V(\varphi)$, there are natural
numbers $n$ and $m$ such that $[q] = n V(\varphi)([1_B]) = n[p]$  and
$[p-q] = m V(\varphi)([1_B]) = m[p]$. In particular, $[p] \le [q]$ and
$[p] \le [p-q]$, and so $2[p] \le [q]+[p-q] = [p]$, which entails that
$p$ is properly infinite.
\end{proof}

\noindent We shall also need the following
lemmas for the proof 
of our main results on divisibility in \Cs s of real rank zero. 

\begin{lem} \label{lm:div1}
Let $M$ be a conical refinement monoid (or a monoid with the Riesz
decomposition property), let $n$ be a
natural number, and let 
$x,y$ in $M$ be such that $y \le x $, $x \propto y$, and $M$ is weakly
divisible of degree $n$ at $y$. It follows that $M$ is weakly
divisible of degree $n$ at $x$.
\end{lem} 

\begin{proof} Find $u \in M$ such that $x = y+u$, and find natural
  numbers $r, n_1, \dots, n_r$, with $n_j \ge n$ for all $j$, and
  elements $y_1, \dots, y_r$ in $M$ 
  such that $y = n_1y_1 + \cdots + n_ry_r$. The assumption $x \propto
  y$ implies that $u \propto y_1 + \cdots + y_r$, and so there is a
  natural number $k$ such that $u \le k(y_1 + \cdots + y_r)$. 

It suffices to show that whenever 
$k$ is a natural number, 
$x,u,y,y_1,$$ \dots, y_r \in M$, and $n_1, \dots, n_r \in \mathbb{N}$
are such that 
$n_j \ge n$ for all $j$, $x =  y + u$, $y= n_1y_1 + \cdots + n_ry_r$,
and $u \le k(y_1 + \cdots + y_r)$, 
then there is a natural number $s$ and there are $u',y',y'_1, \dots, y'_s
\in M$ and $n'_1, \dots, n'_s \in \mathbb{N}$  such that $n_i' \ge n$
for all $i$, 
$x = y' + u'$, $y'= n'_1y'_1 + \cdots + n'_sy'_s$, and $u' \le (k-1)(y'_1 +
\cdots + y'_r)$. The proof of the lemma is completed after $k$
such reductions.

Use the refinement (or the decomposition) property on the inequality
$u \le k(y_1 + \cdots 
+ y_r)$ to find elements $y_{ji} \in M$ with 
$$u = \sum_{j=1}^r \sum_{i=1}^k y_{ji}, \qquad  y_{ji} \le y_j, \; \;
j=1, \dots, r, \, i = 1, \dots, k.$$
Find $z_j \in M$ with $y_{j1}+z_j = y_j$ for $j = 1, \dots, r$. 
Put $s = 2r$, and put
\begin{eqnarray*}
y' & \eqdef & y + y_{11} + y_{21} + \cdots + y_{r1} \\
&=& n_1(y_{11} + z_1) + n_2(y_{21} + z_2) + \cdots +  n_r(y_{r1} +
z_r) +  y_{11} + y_{21} + \cdots + y_{r1} \\
& = & (n_1+1)y_{11} + (n_2+1)y_{21} + \cdots + (n_r+1)y_{r1} + n_1 z_1
+ n_2 z_2 + \cdots + n_r z_r \\
&=& n'_1y'_1 + n'_2y'_2 + \cdots + n'_sy'_s, 
\end{eqnarray*}
where 
$$y_j' = \begin{cases} y_{j1}, & j = 1, \dots, r,\\ z_{j-r}, & j=r+1,
  \dots, s, \end{cases} \qquad 
n_j' = \begin{cases} n_j+1, & j = 1, \dots, r,\\ n_{j-r}, & j=r+1,
  \dots, s. \end{cases}
$$
Moreover,  
$$u' \; \eqdef \;  \sum_{j=1}^r \sum_{i=2}^k y_{ji} \; \le \; \sum_{j=1}^r
\sum_{i=2}^k y_j \; = \; (k-1) \sum_{j=1}^r y_j \; = \;
(k-1)\sum_{j=1}^s y'_j,$$ 
and $x = y'+u'$ as desired. 
\end{proof} 

\noindent We shall a couple of times use the following lemma, whose
proof is verbatim identical to the proof of \cite[Lemma~9.8]{eff}
by Effros, and which therefore is omitted. In the formulation of Effros'
lemma, it is required that 
the ideal $I$ below is an AF-algebra, but inspection of the proof
shows that we only need $I$ to be of real rank zero. 

\begin{lem} \label{lm:div2}
Let $A$ be a \Cs{}, let $I$ be a closed two-sided ideal in $A$, such
that $I$ is of real rank zero, and let $\pi \colon A \to A/I$ be the quotient
mapping. Then each \sh{} $\mu \colon B \to A/I$, where $B$ is a finite
dimensional \Cs{}, lifts to a \sh{} $\lambda \colon B \to A$, i.e.
$\mu = \pi \circ \lambda$. 
\end{lem}

\begin{lem} \label{lm:div3} 
Let $A$ be a unital \Cs{} of real rank zero, and let $n$ be a natural
number. Let $\cP_n$ be the set of projections $p \in A$ for
which $A$ is weakly divisible of degree $n$ at $p$. Suppose $\cP_n$
is full in $A$ (i.e., is not contained in any proper two-sided ideal
in $A$). Then $A$ is weakly divisible of degree $n$ at its unit $1_A$.
\end{lem}

\begin{proof}  We show that $\cP_n$ contains a full projection. 

The assumptions imply that $1_A$ belongs to the
closed---and hence to the algebraic---ideal in $A$ generated by
$\cP_n$. It follows that $1_A$ belongs to the (closed) two-sided ideal
generated by some finite set  of projections in
$\cP_n$.

To prove that $\cP_n$ contains a full projection, it suffices to show
that if $k \ge 2$ and 
$p_1, \dots, p_k$ are projections in $\cP_n$, then there are
projections $q_1, \dots, q_{k-1}$ in $\cP_n$, such that the two sets
$\{p_1, \dots, p_k\}$ and $\{q_1, \dots, q_{k-1}\}$ generate the same
closed two-sided ideal in $A$. 

Let $I$ be the closed two-sided ideal in $A$ generated by
$p_{k-1}$. If $p_k$ belongs to $I$, then we can take $q_j=p_j$ for
$j=1, \dots, k-1$. Suppose that $p_k$ does not belong to $I$. Let $\pi
\colon A \to A/I$ be the quotient mapping. Then $A/I$ is weakly divisible of
degree $n$ at $\pi(p_k)$ (and $\pi(p_k) \ne 0$). Hence there
is a finite 
dimensional \Cs{} $B = M_{n_1} \oplus \cdots \oplus M_{n_r}$, with
$n_j \ge n$ for all $j$, and a unital \sh{} $\mu \colon B \to
\pi(p_k)(A/I)\pi(p_k)$. Let $\pi' \colon (1-p_{k-1})A(1-p_{k-1}) \to
A/I$ be the restriction of $\pi$, and notice that $\pi'$ is
surjective. It follows from Lemma~\ref{lm:div2} that $\mu$ lifts to a
\sh{} $\lambda \colon B \to   (1-p_{k-1})A(1-p_{k-1})$. The two
projections $p_{k-1}$ and $\lambda(1_B)$ are mutually orthogonal and
both belong to $\cP_n$,
so the projection
$q_{k-1}=p_{k-1}+\lambda(1_B)$ 
belongs to $\cP_n$. Let $J$ and $J'$ be the closed two-sided ideals in
$A$ generated by $\{p_{k-1},p_k\}$
and $q_{k-1}$, respectively.
Then $p_{k-1}$ belongs to $J'$ (and clearly also to
  $J$), so $I$ is contained in both $J$ and $J'$. The quotients $J/I$
  and $J'/I$ are generated by $\pi(p_k)$
and $\pi(q_{k-1})$, respectively; but 
$$\pi(q_{k-1}) = (\pi \circ \lambda)(1_B) = \mu(1_B) = \pi(p_k),$$
and this proves that $J/I = J'/I$, whence $J=J'$.
We conclude that $\{p_1, \dots, p_{k-2},q_{k-1}\}$ generates the same
closed two-sided ideal in $A$ as $\{p_1, \dots, p_{k}\}$.

We now have a full projection $p \in \cP_n$. Phrased in the language
of monoids, $[p] \le [1_A] \propto [p]$, and $V(A)$ is weakly divisible of
degree $n$ at $[p]$. Lemma~\ref{lm:div1} implies that $V(A)$ is
weakly divisible of degree $n$ at $[1_A]$, and this in turn implies that $A$
is weakly divisible of degree $n$ at $1_A$.
\end{proof}

\begin{prop} \label{prop:divisible2}
Let $A$ be a separable $C^*$-algebra of real rank zero,
 let $p$ be a projection in $A$, and let $n$ be a natural number. Then $A$ is
weakly divisible of degree $n$ at $p$
if and only if the corner $pAp$ has no (non-zero) representation of
dimension less than $n$.
\end{prop}

\begin{proof} Suppose first that $pAp$ has a non-zero (possibly non-faithful)
  representation $\pi$ on a 
  Hilbert space of dimension $m$, and suppose that
there is a unital \sh{} from $B=M_{n_1} \oplus \cdots \oplus M_{n_r}$
into $pAp$, where $n_j \ge n$ for all $j$. 
The representation $\pi$ will then restrict to a non-zero representation of
$B$, but this is possible only when $m \ge \min\{n_1,\dots,n_r\}$. 
This proves the ``only if'' part of the proposition.

Suppose now that $pAp$ has no (non-zero) representation of dimension
less than $n$. For ease of notation, and upon replacing $A$ with $pAp$, we
can assume that $A$ is unital and that $p=1_A$. Let $\cP_n$ be the set
of all projections $q$ 
in $A$ for which $A$ is weakly divisible of degree $n$ at $q$. Let $I$
be the closed two-sided ideal in $A$ generated by $\cP_n$. We claim
that $I=A$, and this will complete the proof by Lemma~\ref{lm:div3}.

Suppose, to reach a contradiction, that $I \ne A$. Then $I$ is
contained in a maximal ideal $J$ of $A$. Let $\pi \colon A \to A/J$ be
the quotient mapping. The quotient $A/J$ is
separable, simple, unital, and of real rank zero. If $A/J$ is finite
dimensional, then it is isomorphic to $M_m$ for some $m \ge n$ (by the
assumption that $A$ has no representation of dimension less than
$n$). If $A/J$ is infinite dimensional, then it is weakly divisible by
Proposition~\ref{prop:div1}. In either case there is a finite
dimensional \Cs{} $B= M_{n_1} \oplus \cdots \oplus M_{n_r}$, with $n_j
\ge n$ for all $j$, and a unital \sh{}
$\mu \colon B \to A/J$. By Lemma~\ref{lm:div2},
$\mu$ lifts to a \sh{} $\lambda \colon B \to A$. This entails that $q=
\lambda(1_B)$ belongs to $\cP_n$ and hence to $I \subseteq J$, thus
yielding the contradiction
$0 = \pi(q) = (\pi \circ \lambda)(1_B) = \mu(1_B) =
1_{A/J} \ne 0$.
\end{proof}

\noindent If $A$ is weakly divisible of degree $n$ at a projection
$p$, then there
is a full (but not necessarily unital) $^*$-homomorphism from $M_n$
into $pAp$. (A $^*$-homomorphism is called \emph{full} if the closed
two-sided ideal generated by its image is the entire $C^*$-algebra.)
By Proposition~\ref{prop:divisible2} such a $^*$-homomorphism
exists whenever $pAp$ has no representation of dimension less than
$n$, provided that $A$ is of real rank zero. In particular, there is a
full $^*$-homomorphism $M_2 \to pAp$ precisely when $pAp$ has no
character.

In the non-real rank zero case, there are
infinite dimensional simple unital $C^*$-algebras with no projections
other than $0$ and $1$, and there is no (full) embedding of $M_2$ into
such a $C^*$-algebra.

It is not known for
which unital, non-real rank zero algebras $A$ there exists a full
$^*$-homomorphism from $M_n \otimes C_0((0,1])$ into $A$ (even for
$n=2$). It is known,
however, that absence of characters of $A$ is not a sufficient condition to
ensure the existence of a full $^*$-homomorphism from $M_2 \otimes
C_0((0,1])$ into $A$, but absence of
finite dimensional representations could be sufficient. This problem
is referred to as the \emph{Global Glimm problem}, and it has been
considered in the study of
purely infinite $C^*$-algebras, see for example \cite{KirRor} and
\cite{BlaKir}.

\newpage
\begin{thm} \label{thm:divisibility}
Let $A$ be a separable $C^*$-algebra of real rank zero. Then the
following are equivalent:
\begin{enumerate}
\item $A$ is weakly divisible (cf.\ Definition~\ref{def:1}).
\item No non-zero corner $pAp$ (where $p$ is a projection on $A$)
  admits a character.
\item There is no representation $\pi$ of $A$ on a Hilbert space $H$
  for which $\pi(A) \cap {\mathcal{K}}(H) \ne \{0\}$ (where
  ${\mathcal{K}}(H)$ denotes the algebra of compact operators on $H$).
\end{enumerate}
\end{thm}

\begin{proof} (i) $\Rightarrow$ (iii). Suppose that $\pi$ is a
  representation of $A$ on a Hilbert space $H$ such that $\pi(A) \cap
  {\mathcal{K}}(H) \ne \{0\}$. Let $I$ be the kernel of $\pi$, and put
  $J = \pi^{-1}({\mathcal{K}}(H))$. Then $I$ and $J$ are closed
  two-sided ideals, $I \subset J$, and $J/I$ is isomorphic to  $\pi(A)
  \cap {\mathcal{K}}(H)$. The property of having real rank zero passes
  to ideals and quotients, so $J/I$, and hence  $\pi(A) \cap
  {\mathcal{K}}(H)$ are of real rank zero. Take a non-zero projection
  $p$ in in  $\pi(A) \cap {\mathcal{K}}(H)$ and use again
  the real rank zero property of $A$ to lift $p$ to a projection $q$
  in $A$. There is no unital $^*$-monomorphism $M_n \oplus M_{n+1} \to
  p\pi(A)p = \pi(qAq)$ for $n > \dim(p)$, and hence there is no
  unital $^*$-monomorphism $M_n \oplus M_{n+1} \to qAq$ for those $n$.
  This shows that $A$ is not weakly divisible at $q$, and therefore $A$ is
  not weakly divisible.

(iii) $\Rightarrow$ (ii). Suppose that $p$ is a non-zero projection in
$A$ and that $\rho \colon pAp \to \C$ is a character. Put
$\overline{\rho}(a) = \rho(pap)$ for $a$ in $A$. Then
$\overline{\rho}$ is an extremal state on $A$ which extends
$\rho$. Let $(\pi, H,  \xi)$ be the GNS representation of $A$ that
arises from the state $\overline{\rho}$. Then $\pi$ is irreducible,
$\pi(p)\pi(A)\pi(p) = \pi(pAp)$ is abelian, and $\pi(p) \ne 0$
(the latter because $\langle \pi(p) \xi, \xi \rangle =
\overline{\rho}(p)=1$). This entails that $\pi(p)$ is a 1-dimensional
projection, and therefore $\pi(A) \cap {\mathcal{K}}(H) \ne \{0\}$.

(ii) $\Rightarrow$ (i). This follows from
Proposition~\ref{prop:divisible2} (with $n=2$) together with
Lemma~\ref{divisible1}.
\end{proof}

\noindent Condition (iii) is equivalent to the statement that there is
no pair of closed two-sided ideals $I \subset J \subseteq A$ such that
$J/I$ is isomorphic to a sub-$C^*$-algebra of the compact operators
$\mathcal{K}$.

\end{document}